\documentclass[11pt,reqno]{amsart}
\usepackage[margin=1.0in]{geometry}
\usepackage{bm,bbm,amssymb,xcolor,mathtools}

\usepackage[colorlinks=true]{hyperref} %
\hypersetup{citecolor = blue}
\usepackage[noabbrev,capitalize]{cleveref}

\usepackage{bbm}
\usepackage{cite}
\usepackage{psfrag}
\usepackage{cite}
\usepackage{color,soul}
\usepackage{breakcites}     %
\usepackage{bigints}
\usepackage{algorithm}
\usepackage[noend]{algpseudocode}

\usepackage{subcaption}
\usepackage{enumitem}
\setlist{leftmargin=1.6em}

\newtheorem{example}{Example}
\newtheorem{assumption}{Assumption}
\newtheorem{remark}{Remark}
\newtheorem{theorem}{Theorem}
\newtheorem{proposition}{Proposition}
\newtheorem{lemma}{Lemma} 
\newtheorem{corollary}{Corollary}
\usepackage[draft]{graphicx}
\usepackage{xcolor}
\usepackage{svg}
\usepackage{algorithm}

\makeatletter
\def\subsubsection{\@startsection{subsubsection}{3}%
  \z@{.5\linespacing\@plus.7\linespacing}{-.5em}%
  {\normalfont\bfseries}}
\makeatother

\newcommand{\sH}[0]{\mathsf{H}}

\newcommand{\supp}{\mathrm{spt}}

\renewcommand{\vec}[0]{\text{vec}}

\newcommand{\vertiii}[1]{{\left\vert\kern-0.25ex\left\vert\kern-0.25ex\left\vert #1 
    \right\vert\kern-0.25ex\right\vert\kern-0.25ex\right\vert}}

\newcommand{\cC}{\mathcal{C}}

\newcommand{\cI}{\mathcal{I}}

\newcommand{\cK}{\mathcal{K}}

\newcommand{\cO}{\mathcal{O}}

\newcommand{\RR}{\mathbb{R}}

\newcommand{\eqrefprimal}{(\hyperref[eq:entropicQP]{$\mathrm{QP}_0$}) }
\newcommand{\eqrefLP}{(\hyperref[eq:entropicLP]{$\mathrm{LP}_0$}) }

\DeclareMathOperator{\argmin}{argmin}
\DeclareMathOperator{\argmax}{argmax}

\DeclareMathOperator{\dist}{dist}
\DeclareMathOperator{\conv}{conv}

\DeclareMathOperator{\range}{range}
\DeclareMathOperator{\Int}{int}
\DeclareMathOperator{\bd}{bd}
\DeclareMathOperator{\op}{op}

\newcommand{\rM}{\mathrm{M}}
\newcommand{\rA}{\mathrm{A}}

\newcommand{\rP}{\mathrm{P}}
\newcommand{\rD}{\mathrm{D}}
\newcommand{\rB}{\mathrm{B}}

\begin{document}
 \title[Approximation Analysis of the Entropic Penalty in Quadratic Programming]{Approximation Analysis of the\\Entropic Penalty in Quadratic Programming}

\author[V. Karumanchi]{Venkatkrishna Karumanchi}
\address[V. Karumanchi]{School of Operations Research and Information Engineering, Cornell University.}
\email{vk383@cornell.edu}

\author[G. Rioux]{Gabriel Rioux}
\address[G. Rioux]{Department of Mathematics, Imperial College London.}
\email{g.rioux@imperial.ac.uk}

\author[Z. Goldfeld]{Ziv Goldfeld}
\address[Z. Goldfeld]{School of Electrical and Computer Engineering, Cornell University.}
\email{goldfeld@cornell.edu}

  \begin{abstract}
    Quadratic assignment problems are a fundamental class of combinatorial optimization problems which are ubiquitous in applications, yet their exact resolution is NP-hard. To circumvent this impasse, it was proposed to regularize such problems via an  entropic penalty, leading to  computationally tractable proxies. Indeed, this enabled efficient algorithms, notably in the context of Gromov-Wasserstein (GW) problems, but it is unknown how well solutions of the regularized problem approximate those of the original one for small regularization parameters. Treating the broader framework of general quadratic programs (QPs), we establish that the approximation gap decays exponentially quickly for  concave QPs, while the rate for general indefinite or convex QPs can be as slow as linear. Our analysis builds on the study of the entropic penalty in linear programming by leveraging a new representation for concave QPs, which connects them to a family of linear programs with varying costs. Building on these results, we design an algorithm which, given a local solution of the entropic QP, returns a candidate minimizer of the original QP and certifies it. %
    We apply these findings to a general class of discrete GW problems, yielding new variational forms and the first exponentially vanishing entropic approximation bound in the GW literature.

  \end{abstract}
  \date{\today}
  \maketitle
    \section{Introduction}
    \label{sec:introduction}
    
The quadratic assignment problem (QAP) \cite{koopmans1957assignment} is a fundamental combinatorial optimization problem that subsumes  graph matching, the maximum clique problem, and the traveling salesman problem. While QAPs have seen practical applications across economics \cite{heffley1972quadratic,heffley1980decomposition,
koopmans1957assignment}, data analysis \cite{hubert1976quadratic}, and chemistry \cite{forsberg1995analyzing}, among others \cite{loiola2007survey}, their exact resolution is NP-hard \cite{sahni1976p}, in fact, NP-complete in general \cite{Commander2005}. Consequently, several computationally tractable relaxations of the QAP have been proposed in the literature, among which we focus on the so-called softassign QAP \cite{gold1995matching}. The idea is to relax the constraint that the assignment is one-to-one and regularize the objective  with an entropy penalty. %
Specifically, %
following \cite[Equations (2)-(3)]{rangarajan1999convergence}, the softassign QAP can be written as

    \begin{equation}
    \label{eq:softassign}
        \min_{\Sigma\in \mathcal Q_N}\left\{\frac 12 \mathrm{vec}(\Sigma)^{\intercal}\rM\mathrm{vec}(\Sigma)+c^{\intercal}\mathrm{vec}(\Sigma)-\varepsilon \mathsf H(\vec(\Sigma))\right\}, 
        \end{equation}
        where $\mathcal Q_N$ is the set of $N\times N$ bistochastic matrices\footnote{That is, matrices with nonnegative entries whose rows sum to the vector of all $1$'s and similarly for its columns.}, $\mathrm{M}\in\mathbb R^{N^2\times N^2}$ and $c\in\mathbb R^{N^2}$ are specified by the QAP, $\mathrm{vec}( \mathrm L)$ denotes the vectorization of a matrix $\mathrm{L}$, $\varepsilon\geq 0$ is a regularization parameter, and for a vector of length $n$, the entropy term $\mathsf H(x)$ is defined as 
        \[
\mathsf{H}(x)\coloneqq
\begin{cases}
    -\sum_{i=1}^nx_i\log(x_i), &\text{if }x\geq 0,
    \\
    -\infty, &\text{otherwise},
\end{cases}
\]
with the convention that $0\log(0)=0$. Efficient algorithms for solving \eqref{eq:softassign} were proposed and analyzed in \cite{gold1995softassign,rangarajan1996novel,rangarajan1999convergence}. These works propose to gradually anneal the regularization parameter with the hope  that the resulting sequence of regularized solutions provides a good approximation of a solution to the original problem in the limit. Similar regularization techniques were also proposed in the study of Gromov-Wasserstein (GW) problems, which are intimately related to QAPs, again leading to tractable algorithms for the entropic proxy  \cite{peyre2016gromov,solomon2016entropic}. Despite the popularity of the entropic penalty for computation, no theoretical results quantifying the rate at which the regularized solutions converge are currently available. %

\subsection{Contributions}

This work addresses this gap by studying the general question of how solutions  of an arbitrary entropically penalized quadratic program (QP),
    \begin{equation}
        \label{eq:entropicQP}
        \tag{$\mathrm{QP}_{\varepsilon}$}
        \begin{aligned} 
        \text{minimize}\quad &\frac{1}{2} x^{\intercal} \rM x +c^{\intercal} x - \varepsilon \mathsf{H}(x)  
        \\
        \text{subject to}\quad &\rA x= b,
         \\&\phantom{\rA}x\geq 0,
        \end{aligned}
    \end{equation}
    approximate those of its unregularized counterpart (\hyperref[eq:entropicQP]{$\mathrm{QP}_0$}), obtained by setting $\varepsilon=0$ above, in the limit of vanishing regularization. We make precise that $\rM\in\mathbb R^{n\times n}$ is a symmetric matrix\footnote{Symmetry can be assumed without loss of generality as $x^{\intercal}(\rM+\rM^{\intercal})x = 2x^{\intercal}\rM x$.} and work under the following mild assumption throughout.
 \begin{assumption} 
 \label{assn}
The feasible set $\mathcal K\coloneqq \{ x\in\mathbb R^n:\rA x= b,x\geq 0\}$ is nonempty and bounded. Moreover,  \eqref{eq:entropicQP} is nontrivial in the sense that $\frac{1}{2} x^{\intercal}\rM x +c^{\intercal} x$ is nonconstant on $\mathcal K$.
\end{assumption}

Our first main result concerns the rate at which convergence of solutions occurs. 

\begin{theorem}[QP exponential approximation rate; informal] 
\label{thm:informalRate}
For a negative semidefinite matrix $\rM\in\mathbb R^{n\times n}$, if $x_\varepsilon^{\star}$ solves \eqref{eq:entropicQP}, there exists $C_1,C_2,C_3>0$ which depend on $\cK, \rM,$ and $c$, such~that 
    \[
    \dist\left(x_\varepsilon^{\star}, \argmin_{x \in \cK}\left\{\frac{1}{2}x^\intercal \rM x + c^\intercal x\right\}\right) \leq C_1\exp(-{C_2}/{\varepsilon})\text{ for every }0<\varepsilon\leq C_3.
    \]
\end{theorem}
A precise statement, including bounds on each constant, is provided in \cref{thm:EntropicQP}. As a direct consequence of this result, we obtain that $\frac 12 (x_\varepsilon^{\star})^{\intercal}\rM x_\varepsilon^{\star} + c^{\intercal}x_\varepsilon^{\star}$ %
converges towards the optimal value of \eqrefprimal at the same exponential rate, see \cref{cor:cost}. Interestingly, the exponential rate in \cref{thm:informalRate} is not guaranteed to hold unless the QP is concave. Indeed, in \cref{sec:slowRates}, we furnish examples of convex and indefinite QPs for which the regularized solutions converge to their unregularized counterparts at the rate of $\Omega(\varepsilon)$. 

An important component of our analysis is to show that, when $\rM =-\rB^{\intercal}\rB$ for some $\rB\in\mathbb R^{r\times n}$, \eqref{eq:entropicQP} is equivalent to the optimization problem 
\begin{equation}
\label{eq:variationalIntro}
\inf_{u\in\mathbb R^r}g_{\varepsilon}(u),\quad g_{\varepsilon}(u)\coloneqq  \frac{1}{2}\|u\|^2+\min_{x\in\mathcal K} \left\{\left( c-\rB^{\intercal}u\right)^{\intercal}x-\varepsilon \mathsf{H}(x)\right\},
\end{equation}
for every $\varepsilon\geq 0$
in the sense that solutions to one problem can be constructed from solutions of the other. The expression \eqref{eq:variationalIntro} trades off the minimization of a quadratic form over $\mathcal K$ with the unconstrained minimization of the sum of the squared norm with the optimal value function of a linear program (LP) with feasible set $\mathcal K$ and cost dependent on $u$. This representation effectively follows from the fact that the convex quadratic form $-\frac{1}{2}x^{\intercal}\rM x-c^{\intercal}x$ coincides with its biconjugate. %
The remainder of the proof of \cref{thm:EntropicQP} leverages the exponential rate of convergence of solutions of entropic LPs towards their unregularized counterparts, see \cref{thm:entropicRateLP} for a precise statement. Crucially, as the cost vectors in the LP in \eqref{eq:variationalIntro} vary, \cref{thm:entropicRateLP} does not directly apply, as the rate in that result may degenerate when considering sequences of cost vectors. To rectify this, we derive analytic properties of  solutions of \eqref{eq:variationalIntro}  that guarantee that the LP rates do not degenerate.

In \cref{sec:cert}, we use the variational representation \eqref{eq:variationalIntro} and the exponential rate of convergence of regularized solutions to develop an algorithm which maps approximate  critical points of $g_{\varepsilon}$ to candidate local minimizers of $g_0$  and certifies if the candidate is indeed locally minimal. As noted previously, efficient algorithms for approximately solving the softassign QAP and the entropically regularized GW (EGW) problem are available. When the resulting problems are concave,\footnote{\label{foot:softassign}Examples of concave QPs, regularized or not, are abundant. For instance, the softassign QAP modifies the QAP objective so that the quadratic form in \eqref{eq:softassign} is concave. Intuitively, this should ensure that solutions of \eqref{eq:softassign} recover a feasible point for the QAP as $\varepsilon\downarrow 0$, see Section 2 in \cite{rangarajan1999convergence}. Also, the quadratic GW problem is always concave, see \cref{rem:2-2}.} our method enables lifting these entropic approximations to a candidate local minimizer for the original problem and to certify if the candidate is or is not locally minimal.

In \cref{sec:GW}, we apply our general theory to the problem of GW alignment, which has driven recent progress in machine learning
\cite{alvarez2018gromov, bunne2019learning,sejourne2021unbalanced,yan2018semi}, single-cell genomics \cite{blumberg2020mrec,cao2022manifold,
demetci2020gromov}, and object matching \cite{chen2020graph,koehl2023computing,memoli2009spectral,petric2019got,xu2019scalable,xu2019gromov}. The broad applicability of GW problems has motivated recent theoretical developments on their  computational \cite{peyre2016gromov,scetbon2022linear,solomon2016entropic,rioux2023entropic} and statistical \cite{zhang2024gromov,rioux2024limit} properties, in which the effect of entropic regularization is of a particular interest. Despite these advances, the fundamental question of how well solutions of EGW problems approximate solutions of the GW problem remains largely unanswered. To our knowledge, the only result in this direction is Proposition 1 in \cite{zhang2024gromov}, which asserts that the entropic approximation gap for the Euclidean GW problem with quadratic cost vanishes at the generic rate of $O\left(\varepsilon\log(1/\varepsilon)\right)$. While this bound holds in great generality, the rate is governed by the speed at which the entropy term vanishes, and thus cannot leverage specific structure of the underlying distributions (e.g., finitely supported, corresponding to the setting herein). By contrast, the quantitative approximation rates for concave QPs we derive yield exponential approximation rates for both the optimizer and value function of a broad class of discrete GW problems, far beyond the quadratic setting, see  \cref{cor:costGW}. %
To better contextualize these results, we provide a brief overview of optimal transport (OT) and GW problems in  \cref{sec:primerOT}.

\subsection{Literature review}
\label{sec:literatureReview}

As aforementioned, the softassign QAP serves as an important  example of an entropically penalized QP. The interest in softassign QAPs stems from their application to solving hard combinatorial problems including graph matching, the traveling salesman problem, and graph partitioning  via efficient approximation methods as considered in
\cite{gold1995matching,gold1995softassign,gold2002graduated,rangarajan1996novel,rangarajan1999convergence}. Beyond QAPs, entropically regularized QPs  have found applications to transportation planning problems \cite{fang1995linearly} with the aim of improving upon previous approaches based on solving entropically regularized LPs \cite{brice1989derivation,erlander1990efficient,tomlin1971mathematical,tomlin1968traffic}. Numerical optimization of entropic QPs was also studied in \cite{fang1993unconstrained,preda2009convex} as a means to solve convex QPs using unconstrained optimization techniques.

The developments in the present work are inspired by the literature on entropic approximation of LPs. Concretely, 
    that line of work studies how well the entropically regularized LP 
    \begin{equation}
\label{eq:entropicLP}
\tag{$\mathrm{LP}_{\varepsilon}$}
\begin{aligned}
\text{minimize}\quad &c^{\intercal} x{-\varepsilon\mathsf{H}(x)}  
        \\
        \text{subject to}\quad &{\rA} x = b.
        \\&\phantom{\rA}x\geq 0.
        \end{aligned}
\end{equation}    
approximates its unregularized counterpart \eqrefLP\hspace{-.25em}. The added entropy term simultaneously enforces the constraint $x \geq 0$ and serves as a strictly convex regularizer. Leveraging this convexity, \cite{fang1993linear} proposed an unconstrained  approach to solving \eqref{eq:entropicLP} with a quadratic rate of convergence, utilizing Newton's method and the steepest descent method on the geometric dual problem. Another important application of entropic LPs  arises in OT theory, where \cite{cuturi2013sinkhorn} proposed Sinkhorn's diagonal scaling algorithm \cite{sinkhorn1967diagonal} to solve the entropic OT (EOT) problem. For transportation between distributions supported on $N$ points, this approach admits a computational complexity of $O(N^2)$ (cf. e.g., Section 4.3 in \cite{peyre2019computational}), which is a notable speedup over the $O(N^3)$ complexity of solving the standard OT problem via the network simplex method \cite{orlin1997polynomial,tarjan1997dynamic} (see also Section 3.5 in \cite{peyre2019computational}). To justify the entropic framework as a valid proxy, \cite{cominetti1994asymptotic} showed that the solutions $(x^\star_\varepsilon)_{\varepsilon >0}$ of \eqref{eq:entropicLP} converge to a particular solution of \eqrefLP exponentially quickly, in the asymptotic sense, as $\varepsilon\downarrow 0$. It was not until \cite{weed2018explicit} that a quantitative version of this rate was derived, see \cref{sec:entropicLP} ahead for details.

The study of EGW problems has followed a similar trajectory. Namely, it was proposed to entropically regularize GW problems as a path towards efficient algorithms \cite{peyre2016gromov,scetbon2022linear,solomon2016entropic}. While the added entropy no longer convexifies the original problem, these works propose to linearize the EGW problem and apply iterative methods to solve it so that each iteration only requires the resolution of an EOT problem. Applying the methods from \cite{peyre2016gromov,solomon2016entropic} to align two distributions on $N$ points bears a per iteration complexity of $O(N^3)$ stemming from the cost of updating the cost matrix. Under certain low-rank assumptions on the cost matrix, this computational cost can be reduced to $O(N^2)$, as illustrated in \cite{scetbon2022linear}. However, none of these works derived non-asymptotic convergence rates for the proposed algorithms. An alternative approach for the Euclidean quadratic and inner product EGW problem was provided in \cite{rioux2023entropic} with a per iteration complexity of $O(N^2)$ along with non-asymptotic convergence guarantees. This algorithm leverages a  variational form for those GW problems as derived in \cite{zhang2024gromov}; the equivalent formulation for concave QPs \eqref{eq:variationalIntro} is inspired by this earlier work.

\section{Notation and Preliminaries}
   
    \subsection{Notation}
    For a point $z\in\mathbb R^d$, its Euclidean norm is denoted by $\|z\|$ and $\|z\|_1$ is its $1$-norm. %
      Given a set $\cC\subset\mathbb R^d$, we write $\Int(\cC)$, $\bd(\cC)$, and $\conv(\cC)$ for its interior, boundary, and convex hull, respectively. If $\cC$ is a finite set, we denote its cardinality by $|\cC|$.     
      The point-to-set distance between $z\in\mathbb R^d$ and $\cC$ is given by $\dist(z,\cC)=\inf_{z'\in \cC}\|z-z'\|.$ 
The indicator function of $\cC$ is
    \[
        \delta_{\cC}:x\in\mathbb R^d \mapsto \begin{cases}
    0 & \text{if } x \in \cC, \\
    +\infty & \text{otherwise.} 
    \end{cases}  
    \]

      Given a matrix $\mathrm{A}\in\RR^{m \times n}$, its operator norm is given by $\|\mathrm A\|_{\mathrm{op}}=\sup_{\|x\|=1}\|\mathrm Ax\|$.
    For a symmetric positive semidefinite matrix $\rM \in \RR^{n \times n},$ its pseudoinverse,  $\rM^{\dagger},$ is the unique symmetric positive semidefinite matrix satisfying $\rM\rM^{\dagger}x=\rM^{\dagger}\rM x=x_{\range(\rM)}$, where $x_{\range(\rM)}$ is the orthogonal projection of $x$ onto the range of $\rM$, ${\range(\rM)}.$  

    Throughout the paper, asymptotic notation is always understood as $\varepsilon\downarrow 0$. For functions $f:[0,\infty)\to \mathbb R$ and $g:[0,\infty)\to [0,\infty)$, the notation $f(\varepsilon
    )=O(g(\varepsilon))$ is understood as $|f(\varepsilon)|\leq Cg(\varepsilon)$ for some constant $C$, $f(\varepsilon
    )=o(\varepsilon)$ means that $g(\varepsilon)^{-1}f(\varepsilon)\to 0$ as $\varepsilon\downarrow 0$, and $f(\varepsilon
    )=\Omega(g(\varepsilon))$ indicates that $\limsup_{\varepsilon\downarrow 0
    }g(\varepsilon)^{-1}|f(\varepsilon)|> 0$.

\subsection{Entropic Penalty in Linear Programming}
    \label{sec:entropicLP}
    Recall the entropically penalized LP with regularization strength $\varepsilon\geq 0$, \eqref{eq:entropicLP}, 
    given by
\[
\tag{$\mathrm{LP}_{\varepsilon}$}
\begin{aligned}
\text{minimize}\quad &c^{\intercal} x{-\varepsilon\mathsf{H}(x)}  
        \\
        \text{subject to}\quad &{\rA} x = b,
        \\&\phantom{\rA}x\geq 0.
        \end{aligned}
\]
    As noted previously, the regularized version can be solved more efficiently than the vanilla LP for certain classes of problems. This has motivated the study of the bias introduced when treating \eqref{eq:entropicLP} as a proxy for \eqrefLP\hspace{-.35em}.

    Remarkably, \cite{cominetti1994asymptotic,weed2018explicit} established that, when the feasible set $\mathcal K$ is nonempty and bounded, the unique solutions\footnote{\label{foot:UniquenessEntropic}It is easy to see that the objective is strictly convex and continuous for any $\varepsilon>0$.  Compactness and convexity of $\mathcal K$ then guarantee  existence and uniqueness of solutions to \eqref{eq:entropicLP}.} $x_{\varepsilon}^{\star}$ to  \eqref{eq:entropicLP} with regularization strength $\varepsilon>0$ converge to the solution set of~\eqrefLP in the point-to-set distance  at an exponential rate as $\varepsilon\downarrow 0$. Precisely, the following quantitative convergence rate  was derived in \cite{weed2018explicit}.\footnote{While the cited reference \cite{weed2018explicit} derived \cref{thm:entropicRateLP} under the point-to-set distance with the $1$-norm, it is more convenient to work with the $2$-norm in the current setting.}

    \begin{theorem}[Exponential convergence rate, Corollary 9 in \cite{weed2018explicit}]
    \label{thm:entropicRateLP} 
    Fix the constants 
   \[
   R_{\sH} \coloneqq \max_{x,x'\in \cK} \left\{\sH(x) - \sH(x')\right\}, \quad R_1 \coloneqq \max_{x \in \cK} \|x\|_{1},
   \] 
    and set $\mathcal O_c\coloneqq \argmin_{x\in\mathcal K}c^{\intercal}x$. Then, if $x^{\star}_{\varepsilon}\in\argmin_{x\in\mathcal K}\left\{c^{\intercal}x-\varepsilon \mathsf H(x)\right\}$, 
   \[
        \mathrm{dist}(x^{\star}_{\varepsilon},\cO_c) \leq 2 R_1 \exp{\left(\frac{-\kappa_c}{\varepsilon R_1} + \frac{R_1 + R_{\sH}}{R_1}\right)} \text{ for every } \varepsilon\leq \frac{\kappa_c}{R_1+R_{\sH}},
   \]
    where 
   $
        \kappa_c\coloneqq \min_{x\in\mathcal V(\mathcal K)\backslash\mathcal O_c}\left\{c^{\intercal}x\right\}-\min_{x\in\mathcal O_c}\left\{c^{\intercal}x\right\}
   $
   and $\mathcal V(\mathcal K)$ is the set of all vertices of $\mathcal K$,
   \end{theorem}

This result serves as an important component of our analysis of the entropic approximation rate for QPs, via a linearization argument that represents the QP in terms of an infimum of a class of LPs with varying costs. To obtain a positive exponential rate for the QP setting we conduct a careful analysis to show that the suboptimality gap $\kappa_c$ does not vanish over the class of LPs.

\subsection{Nonsmooth Optimization} Many standard techniques from optimization theory rely on the assumption that the objective function is sufficiently smooth. In the absence of smoothness, convex analysis provides an elegant toolkit for extending these ideas to the nonsmooth setting. Without smoothness or convexity, significant progress  has been made in the development of nonsmooth analysis techniques for locally Lipschitz continuous functions (i.e., functions which are Lipschitz continuous on every compact set). 

In the convex setting, the analogue of the gradient is given by the  convex subdifferential, which collects all  affine underestimators of the function at a given point. In the case of a locally Lipschitz continuous function,  the Clarke subdifferential \cite{clarke1975generalized} serves as a natural generalization of the gradient which coincides with the convex subdifferential if the function is convex. Given a function $f:X\subset \mathbb R^n\to \mathbb R$, its Clarke subdifferential at $x\in X$  is given by the set 
\[
    \partial f(x)\coloneqq \left\{ \xi \in \mathbb R^n : \limsup_{\substack{y\to x\\t\downarrow 0}}\frac{f(y+tv)-f(y)}{t}\geq \xi^{\intercal}v\text{ for all }v\in X\right\}.
\]
Of note is that  $\partial f(x) = \{\nabla f(x)\}$ if $f$ is continuously differentiable at $x\in X$, illustrating that $\partial f$  properly extends the gradient. Furthermore, if $x$ is locally minimal or maximal for $f$, we have $0\in\partial f(x)$, so that the notion of a stationary point for a smooth function also carries over to this setting. We refer the reader to Chapter 2 of \cite{clarke1990optimization} for a comprehensive treatment of the Clarke subdifferential.

\section{Approximation Gap for Concave Quadratic Programs}
\label{sec:concaveQP}
        The softassign QAP is the entropically regularized QP 
    \begin{equation} 
    \label{eq:softassign2}
    \min_{\Sigma\in \mathcal Q_N}\left\{\frac 12 \mathrm{vec}(\Sigma)^{\intercal}\rM\mathrm{vec}(\Sigma)+c^{\intercal}\mathrm{vec}(\Sigma)-\varepsilon \mathsf H(\vec(\Sigma))\right\},  
        \end{equation} 
     with $\mathcal Q_N$ the set of $N\times N$ bistochastic matrices, $\mathrm{M}\in\mathbb R^{N^2\times N^2}$, and $c\in\mathbb R^{N^2}$ (see \eqref{eq:softassign}). This formulation serves as a relaxation of the standard QAP and has seen applications to approximating solutions of various hard combinatorial problems as discussed in \cref{sec:introduction}. As noted in \cref{foot:softassign}, the matrix $\rM$ is obtained by subtracting $\rho\mathrm{Id}$  from the original QAP cost matrix for some sufficiently large $\rho$  so that $\rM$ becomes negative semidefinite. This modification enforces that solutions of \eqref{eq:softassign2} with $\varepsilon=0$ are feasible for the QAP (see \cite{bazaraa1982use}), while the added entropy term enables efficient algorithms for $\varepsilon>0$. A common belief in this space is that solutions of \eqref{eq:softassign2} with $\varepsilon>0$ approximate their unregularized counterparts well. This section makes this assertion rigorous for general concave QPs, capturing the softassign QAP as a special case.

For a symmetric, positive semidefinite (PSD) matrix $\rM\in\mathbb R^{n\times n}$ and $\varepsilon>0$, consider the entropically regularized concave QP 
\begin{equation}
\label{eq:primalQPND}
    \min_{x\in\mathcal K}\left\{-\frac{1}{2}x^{\intercal}\rM x+c^{\intercal}x-\varepsilon\mathsf{H}(x) \right\}, 
\end{equation}
as a computationally tractable proxy of \eqrefprimal (see \cref{rem:entropic_computation} ahead for a discussion of computational aspects). This section establishes a quantitative relationship between a solution $x_\varepsilon^{\star}$ of \eqref{eq:primalQPND} and the $\argmin$ set of the original (unregularized) QP. To that end, we tie the concave QP to a class of LPs via the following variational representation, which allows leveraging existing bounds on the entropic approximation gap for LPs, see \cref{thm:entropicRateLP}.

    \begin{proposition}[Variational formulation] 
    \label{prop:variationalFormND}
   Suppose, without loss of generality, that $\rM=\rB^{\intercal}\rB$  for some $\rB\in\mathbb R^{r\times n}$ with $r\leq n$.\footnote{Indeed, if $\rM$ is full rank, we may diagonalize it as $\rM=\rP\rD\rP^{\intercal}$ and take $\rB=\rD^{\frac 12}\rP^{\intercal}$. When $r\coloneqq \mathrm{rank}(\rM)<n$, we have $\rM=\tilde {\rP}\tilde{\rD}\tilde{\rP}^{\intercal}$, where $\tilde{\rD}\in\mathbb R^{r\times r}$ is the diagonal matrix containing the top $r$ eigenvalues of $\rM$ and $\tilde{\mathrm{P}}\in \mathbb R^{n\times r}$ has the corresponding eigenvectors in its columns.} Then, for any $\varepsilon \geq 0,$ the optimal value in \eqref{eq:primalQPND} coincides with 
    \begin{equation}
        \label{eq:variationalFormND}
\inf_{u\in\mathbb R^r}\left\{ \frac{1}{2}\|u\|^2+\min_{x\in\mathcal K} \left\{\left( c-\rB^{\intercal}u\right)^{\intercal}x-\varepsilon \mathsf{H}(x)\right\}\right\}.
    \end{equation} Moreover:
   \begin{enumerate}
       \item if $x^{\star}_\varepsilon$ solves \eqref{eq:primalQPND}, then $u^{\star}_\varepsilon=\rB x^{\star}_\varepsilon$ solves \eqref{eq:variationalFormND} and $x^{\star}_\varepsilon\in\argmin_{x\in\mathcal K}\left\{\left(c-\rB^{\intercal}u^{\star}_\varepsilon\right)^{\intercal}x-\varepsilon\mathsf{H}(x)\right\}$,
       \item if $u^{\star}_\varepsilon$ solves \eqref{eq:variationalFormND}, then any $\bar x\in\argmin_{x\in\mathcal K}\left\{\left(c-\rB^{\intercal}u^{\star}_\varepsilon\right)^{\intercal}x-\varepsilon\mathsf{H}(x)\right\}$ solves \eqref{eq:primalQPND}  and  $u^{\star}_\varepsilon=\rB \bar x$.
   \end{enumerate} 
    \end{proposition}

\cref{prop:variationalFormND} is a consequence of the Fenchel-Moreau theorem, which asserts that $f:x\in\mathbb R^n\mapsto \frac{1}{2}x^{\intercal}\rM x-c^{\intercal}x$ coincides with the convex conjugate of its convex conjugate;
see \cref{proof:prop:variationalFormND} for a full derivation.

\begin{remark}[Computational tractability of \eqref{eq:variationalFormND}]\label{rem:entropic_computation}
The variational form in \cref{prop:variationalFormND} offers a pathway for approximately solving concave QPs with and without regularization, which may be of independent interest. This is since it enables recasting the problem as the unconstrained minimization of the objective $\frac 12 \|\cdot\|^2+\min_{x\in\mathcal K}\left\{ (c-\rB^{\intercal}(\cdot))^{\intercal}x-\varepsilon\mathsf{H}(x)\right\}\eqqcolon \frac 12 \|\cdot\|^2+V_{\varepsilon}$, which is the sum of the smooth strongly convex quadratic and the concave function $V_{\varepsilon}$. The proof of \cref{prop:variationalFormND} (see \cref{proof:prop:variationalFormND}) shows that the Clarke subdifferential of the objective is 
\begin{equation}
\label{eq:subdifferential}
   \partial\left(\frac 12 \|\cdot\|^2+ V_{\varepsilon}\right)(u) = u-\rB\argmin_{x\in\mathcal K}\left\{-x^{\intercal}(\rB^{\intercal}u-c)-\varepsilon\mathsf{H}(x)\right\},  
\end{equation}
so that, a priori, $V_{\varepsilon}$ may be nonsmooth. However, as noted in 
\cref{foot:UniquenessEntropic}, for any $\varepsilon>0$, the above $\argmin$ set is a singleton for any $u\in\mathbb R^r$, so that $V_{\varepsilon}$ is continuously differentiable by virtue of the Corollary on p.34 in \cite{clarke1990optimization}. 

This means that the regularized variational form is amenable to optimization via first-order methods. Given the lack of convexity, such algorithms aim to recover critical points of the objective, i.e., $\bar u$ satisfying $\bar u=\mathrm B x_{\bar u}$, where $\{x_{\bar u}\}=\argmin_{x\in\mathcal K}\left\{-x^{\intercal}(\rB^{\intercal}u-c)-\varepsilon\mathsf{H}(x)\right\}$. 
As noted in \cref{sec:literatureReview}, this regularized LP can be solved more efficiently than the unregularized  LP for certain feasible sets $\mathcal K$, such as the set of couplings of two measures.   
These observations were leveraged in \cite{rioux2023entropic} to develop the first algorithms for approximating the entropic quadratic GW distance which are subject to non-asymptotic convergence rate guarantees. 
\end{remark}

\begin{figure}[!tb]
    \centering
    \includegraphics[width=0.7\linewidth]{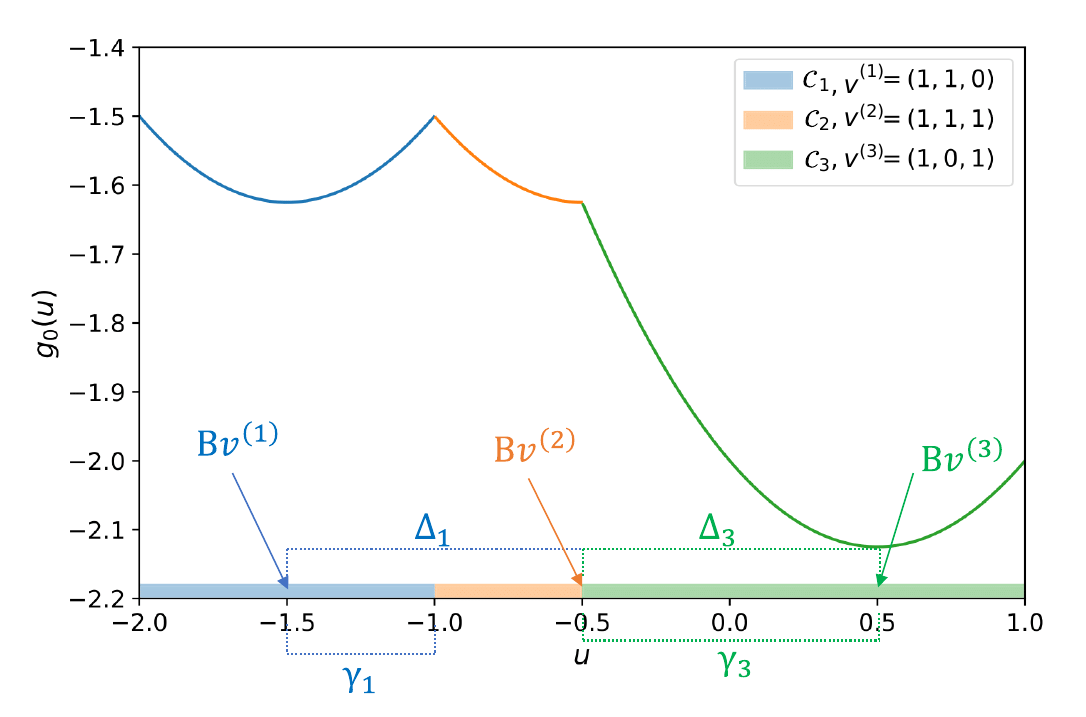}
    \caption{
    \label{fig:variational_form}
    Plot of the objective in the unregularized variational formulation for $\mathrm{M}=zz^{\intercal}$ with $z=(-.5,-1,1)$, $c=(-1,.5,-1)$, and $\mathcal K=[0,1]^3$. The $u$-axis is colored to highlight the regions where a given vertex is optimal. See \cref{thm:EntropicQP} for the definitions of the constants $\Delta_i$ and $\gamma_i.$ }
\end{figure}

While \cref{prop:variationalFormND} provides a connection to LPs, one cannot directly apply \cref{thm:entropicRateLP} since the cost vector $c-\rB^{\intercal}u^{\star}_{\varepsilon}$ featuring in \eqref{eq:variationalFormND} varies with $\varepsilon$. %
To bridge this gap, we study how the optimal vertices of the LP vary as the cost changes. Let $\mathcal V(\mathcal K)=\left\{v^{(i)}\right\}_{i=1}^V$ be the set of vertices of $\mathcal K$ and
\[
    \mathcal C_i\coloneqq \left\{u\in\mathbb R^r:(c-\mathrm{B}^{\intercal}u)^{\intercal} v^{(i)}\leq (c-\mathrm{B}^{\intercal}u)^{\intercal} v^{(j)} \text{ for every $i\neq j$} \right\}
    \]
    denote the set of all $u$'s where the vertex $v^{(i)}$ is optimal for the resulting LP. The set $\{\mathcal C_i\}_{i=1}^V$ is a polyhedral partition of $\mathbb R^r$. Further, let
\[g_\varepsilon:u\in\mathbb R^r\mapsto \frac{1}{2}\|u\|^2 + \min_{x \in \cK}\{(c-\rB^\intercal u )^\intercal x - \varepsilon \sH(x)\}, \text{ for any }\varepsilon\geq 0,\]
be the objective function in the variational form \eqref{eq:variationalFormND}. Note that $g_0$ is piecewise quadratic, as the restriction of $g_0$ to $\mathcal C_i$ is given by $\frac 12 \|\cdot\|^2+(c-\rB^{\intercal}(\cdot))^{\intercal}v^{(i)}$ for each $i\in[V]$, see \cref{fig:variational_form} for an example. The following proposition compiles some properties of $g_0$ that will be useful in the sequel.

\begin{proposition}[On $g_0$ and its minimizers]\label{prop:minimizersG0} Fix $i,j\in[V]$ and assume that there exists $u\in \mathcal C_i\cap \mathcal C_j$. Then either $\mathcal C_i=\mathcal C_j$ or $u\in \mathrm{bd}(\mathcal C_i)\cap \mathrm{bd}(\mathcal C_j)$. Furthermore, if $u\in\mathrm{bd}(\mathcal C_i)$, every neighborhood of $u$ contains a point $u'$ for which $g_0(u')<g_0(u)$. Consequently, every local or global minimizer, $\bar u$, of $g_0$ is such that $\Int(\mathcal C_i)\ni\bar u=\mathrm{B}v^{(i)}$ for some $i\in[V]$.   
\end{proposition}

The proof of the first assertion in \cref{prop:minimizersG0} is a direct consequence of the definition of the regions $\{\mathcal C_i\}_{i=1}^V$. The second assertion follows by characterizing the directional derivative of $g_0$ at a boundary point and showing that there always exists a direction of strict descent. The complete derivation is provided in \cref{proof:prop:minimizersG0}.

\medskip

With these preliminaries in place, we now provide a precise formulation of the informal convergence rate result from \cref{thm:informalRate}.

\begin{theorem}[Approximation gap]\label{thm:EntropicQP}
Let $\cI^{\star}\coloneqq \left\{ 
  i \in [V] : 
  \mathcal C_i \cap \argmin_{u \in \RR^r} g_0(u) \neq \emptyset 
\right\}$, define  
\begin{align*}
\alpha &\coloneqq
\min_{u\in\cup_{i\in [V]\backslash \cI^{\star}}\mathcal C_i}g_0(u)-\min_{u\in\mathbb R^r}g_0(u)>0,\\
\gamma_i&\coloneqq \dist\left(\rB v^{(i)},\mathrm{bd}(\mathcal C_i)\right)>0,\quad i\in \cI^{\star},\\
\Delta_i &\coloneqq \min\Bigg\{\kappa_{c-\rB^{\intercal}\rB v^{(i)}}, \frac{\gamma_i}{2}\min_{\substack{k \in [V] \\ \rB v^{(k)} \neq \rB v^{(i)}}} 
\left\| \rB v^{(k)} - \rB v^{(i)} \right\|\Bigg\}>0,\quad i\in \cI^{\star},
\end{align*}
and set
$ \gamma\coloneqq \min_{i\in \cI^{\star}} \gamma_i$ and $\Delta\coloneqq \min_{i\in \cI^{\star}} \Delta_i$.
Then, for any $0 < \varepsilon < \min\left\{
  \frac{\alpha}{R_{\mathsf H}},\frac{\gamma^2}{8 R_{\mathsf H}},
  \frac{\Delta}{R_1 + R_{\sH}}
\right\}$ and solution $x_{\varepsilon}^{\star}$ of the regularized QP \eqref{eq:primalQPND} with regularization strength \( \varepsilon \), we have that
\[
\dist\left(
  x^\star_\varepsilon,
  \argmin_{x \in \cK}\left\{-\frac{1}{2}x^\intercal \rM x + c^\intercal x\right\}
\right)
\leq
2 R_1 \exp\left(
  \frac{-\Delta}{2\varepsilon R_1}
  + \frac{R_1 + R_{\sH}}{R_1}
\right).
\] 
\end{theorem}

As noted previously, the proof of \cref{thm:EntropicQP}, included in \cref{proof:thm:entropicQP}, leverages the variational form \eqref{eq:variationalFormND} to connect the concave QP to a class of LPs with varying cost. We then demonstrate that, for every $\varepsilon>0$ sufficiently small, $\mathrm{B}x_{\varepsilon}^{\star}\in \mathrm{int}(\mathcal C_i)$ for some $i\in \cI^{\star}$ (where $i$ may depend on  $\varepsilon$) and bound the distance of $\mathrm{B}x_{\varepsilon}^{\star}$ to $\mathrm{bd}(\mathcal C_i)$. As $\mathrm{B}x_{\varepsilon}^{\star}\in \mathrm{int}(\mathcal C_i)$, it holds that $\argmin_{x \in \cK}(c-\rB^\intercal\rB x_{\varepsilon}^{\star})^\intercal x=\argmin_{x \in \cK}(c-\rB^\intercal\rB v^{(i)})^\intercal x$ so that  \cref{thm:entropicRateLP} can be applied to obtain a quantitative bound on the distance of    $x_{\varepsilon}^{\star}$ to $\argmin_{x \in \cK}(c-\rB^\intercal\rB v^{(i)})^\intercal x\subset\argmin_{x\in\mathcal K}\{-\frac 12 x^{\intercal}\rM x+c^{\intercal}x\}$. Finally, we provide a uniform lower bound on the constant $
\kappa_{\left(c-\mathrm{B}^{\intercal}\mathrm{B}x_{\varepsilon}^{\star}\right)}
$ figuring in \cref{thm:entropicRateLP} by leveraging our estimate on $\dist(\rB x_{\varepsilon}^{\star},\mathrm{bd}(\mathcal C_i))$.

\begin{remark}[Improved rate]
   While \cref{thm:EntropicQP} presents the worst-case rate in terms of the constants $\Delta$ and $\gamma$, if $\rB{x}^{\star}_{\varepsilon}\in \mathcal C_i$, the effective rate is governed by $\Delta_i$ and $\gamma_i$. However, we highlight that, much like the constants in \cref{thm:entropicRateLP}, $\Delta_i$ and $\gamma_i$ cannot generally be computed a priori.
\end{remark}

Note that \cref{thm:EntropicQP} is presented in terms of the distance of $x^\star_\varepsilon$ to the set of minimizers of the unregularized QP since it is not guaranteed that a sequence $(x^\star_\varepsilon)_{\varepsilon>0}$ solving \eqref{eq:primalQPND} with regularization $\varepsilon$ is convergent as $\varepsilon \downarrow 0.$ Furthermore, the constant $\Delta_i$ in \cref{thm:EntropicQP} can take the value $+\infty$, but only if every vertex achieves the same for the quadratic form in \eqref{eq:primalQPND}. This cannot occur in the current setting due to \cref{assn}.

As a direct consequence of  \cref{thm:EntropicQP} (namely, by applying the Cauchy-Schwarz inequality), the quadratic form in \eqref{eq:primalQPND} evaluated at  $x^\star_\varepsilon$ converges to its minimum at an exponential rate.
\begin{corollary}\label{cor:cost}
    In the setting of \cref{thm:EntropicQP}, for $0 < \varepsilon < \min\{\frac{\alpha}{2 \beta},\frac{\gamma^2}{8 R_{\sH}},\frac{\Delta}{R_1 + R_{\sH}}\},$ it holds that
    \[0\leq -\frac{1}{2}{(x^\star_\varepsilon)}^{\intercal}\rM x^\star_\varepsilon+c^{\intercal}x^\star_\varepsilon - \min_{x \in \cK}\left\{-\frac{1}{2}x^{\intercal}\rM x+c^{\intercal}x\right\} \leq   K \exp\left(\frac{-\Delta}{\varepsilon R_1}+\frac{R_1 + R_{\sH}}{R_1}\right),\]
    where $K \coloneqq \max\{4R_1^2\|\rM\|_{\op}, 2\|c\|R_1\}.$
\end{corollary}

\section{Certification Algorithm for Regularized Concave Quadratic Programs}\label{sec:cert}

\begin{figure}[!tb]
    \centering
    \includegraphics[width=0.7\linewidth]{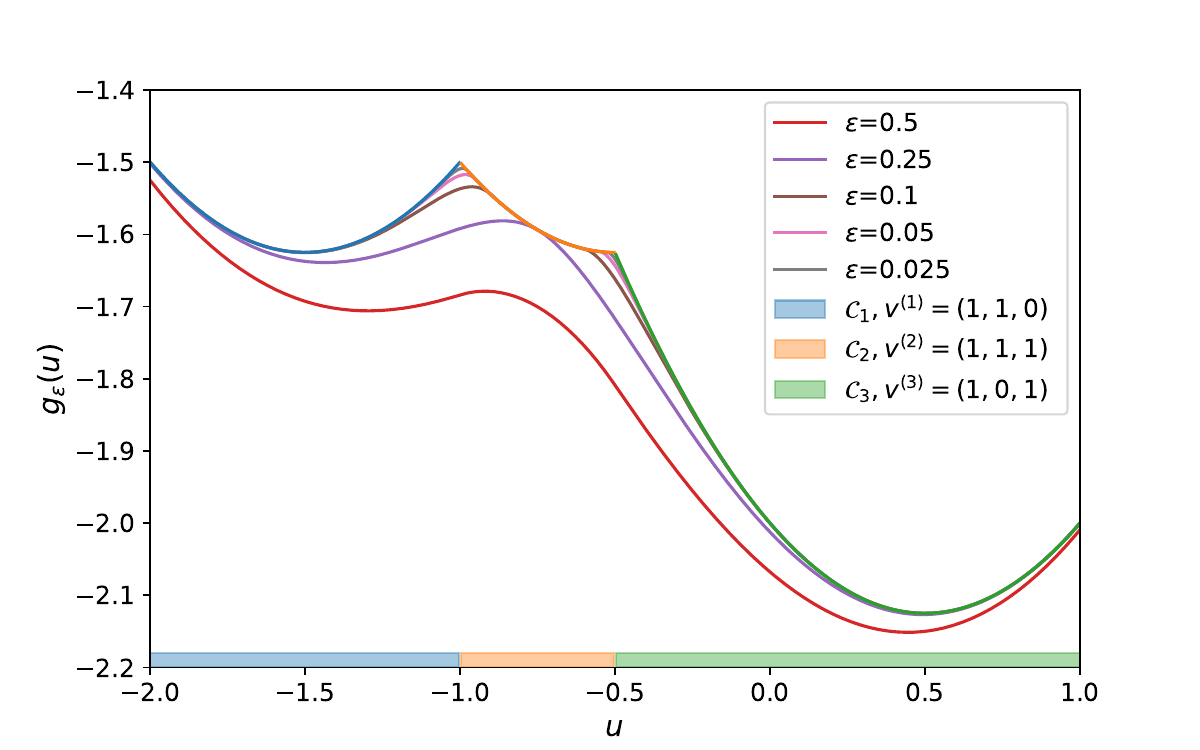}
    \caption{Profiles of $g_{\varepsilon}$ with $\varepsilon\geq 0$ for the example in  \cref{fig:variational_form}. For all  values of $\varepsilon$ considered, $g_{\varepsilon}$ is nonconvex. $g_{\varepsilon}$ is seen to converge quickly to $g_0$ away from the boundary points of the regions $\cC_1,\cC_2$ and $\cC_3,$ and the critical points are seen to concentrate about boundary points of the $\cC_i$ or local minima of $g_0$.}
    \label{fig:figure2}
\end{figure}

While the previous section established convergence of global minimizers of the regularized concave QP, solving this problem for small $\varepsilon>0$ is challenging due to nonconvexity. As an alternative, \cref{rem:entropic_computation} discussed solving the variational problem \eqref{eq:variationalFormND} via first-order methods, since $g_{\varepsilon}$ is smooth for all $\varepsilon>0$. However, gradient-based methods are only guaranteed to return a point $u_\varepsilon\in\RR^r$ with $\|\nabla g_{\varepsilon}(u_\varepsilon)\|\le \delta$ for a prescribed tolerance $\delta$ (a so-called $\delta$-critical point) due to nonconvexity of $g_\varepsilon$; see \cref{fig:figure2} for example.

By studying geometric properties of $g_{\varepsilon}$ for $\varepsilon\geq 0$, 
we develop a procedure that, under suitable conditions, maps $\delta$-critical points of $g_\varepsilon$ to candidate local minimizers of $g_0$. %
Our approach is based on the following result which asserts that once $\varepsilon$ and $\delta$ are sufficiently small, $\delta$-critical points of $g_{\varepsilon}$ are either close to a local minimizer of $g_0$ or to a boundary point of some $\mathcal C_i$ for $i\in[V]$ as illustrated in \cref{fig:figure2}.

\begin{proposition}
\label{prop:nearCriticalPoints}
   For each $i\in[V]$, set $D_i = \dist(\rB v^{(i)},\bd(\cC_i))$ and, for $\gamma>0$, define  
    \[K_i\left(\gamma \right)\coloneqq 
 \min\Bigg\{\gamma \min_{\substack{k\in[V]\\ \rB v^{(k)}\neq\rB v^{(i)} }}\left\|\rB v^{(k)}-\rB v^{(i)}  \right\|,\min_{\substack{k\in[V]\\\rB v^{(k)}=\rB v^{(i)}\\ c^{\intercal}v^{(k)}\neq c^{\intercal}v^{(i)}
    }
    } c^{\intercal}(v^{(k)}-v^{(i)})\Bigg\}.
    \]
Then, for $0<\delta<4\|\rB\|_{\mathrm{op}}R_1$ and $0<\varepsilon< K_i\left( \gamma\right)\left( R_1+R_H-R_1\log\left(\frac{\delta}{4\|\rB\|_{\op}R_1}\right)\right)^{-1},$ we have that 
\begin{enumerate}
    \item if $\rB v^{(i)}\in \Int(\cC_i)$, $0<\gamma<D_i/2$, and %
    $\delta<D_i$,  then each $u\in \cC_i$ with $\|u-\rB v^{(i)}\|\leq ~\frac \delta 2$ is such that $\|\nabla g_{\varepsilon}(u)\|\leq \delta$. Furthermore, a point $u$ satisfies $\|\nabla g_{\varepsilon}(u)\|\leq \delta$ only if (a) $\dist(u,\bd(\cC_i))\geq \gamma$ and $\|u-\rB v^{(i)}\|<\frac{3\delta}{2}$, or (b) $\dist(u,\bd(\cC_i))<\gamma$. Finally, the sets satisfying points (a) and (b) are disjoint once $\delta<D_i/3$.
    \item if $\rB v^{(i)}\not\in\Int(\cC_i)$ and $\delta\leq \frac{2}{3}\gamma$, then $\|\nabla g_{\varepsilon}(u)\|\leq \delta $ only if $\dist(u,\bd(\cC_i))<\gamma$.  
\end{enumerate}
\end{proposition}

The proof of \cref{prop:nearCriticalPoints}, included in \cref{proof:prop:nearCriticalPoints}, follows similar lines to that of \cref{thm:EntropicQP}. Indeed, by leveraging \cref{thm:entropicRateLP}, it is demonstrated that, for small $\varepsilon>0$ and each $u\in\cC_i$ within a prescribed distance from $\bd(\cC_i)$, the solutions of the regularized LP with cost $c-\rB^{\intercal} u$ approximate those of the standard LP with exponential accuracy. This readily implies that  $\|\nabla g_{\varepsilon}(u)-\nabla g_0(u)\|$ can be made arbitrarily small uniformly over all such $u$ by choosing $\varepsilon>0$ small enough.

Item (1) of \cref{prop:nearCriticalPoints} implies that certain $\delta$-critical points of $g_{\varepsilon}$ serve as good proxies for local minimizers of $g_0$. However, some sequences of local minimizers of $g_{\varepsilon}$ may fail to converge to a local minimizer of $g_0$, but rather  converge to saddle points. We illustrate this potential pitfall using the following example, which motivates the upcoming certification algorithm. %

\noindent
\begin{minipage}{.64 \textwidth}
\begin{example}
\label{ex:convergencecritical}

    Consider the one-dimensional concave quadratic program $\min_{x\in [1,2]}\left\{-\frac 12 x^2+x\right\}$. The corresponding regularized variational problem is given by $\min_{u\in\mathbb R} g_{\varepsilon}(u)$ with
    \[
    g_{\varepsilon}(u)=\frac{1}2 u^2 +\min_{x\in[1,2]}\left\{ (1-u)x-\varepsilon \mathsf H(x)\right\}.
    \]
    Solving the minimization over $x\in[1,2]$, the objective can be expressed in closed form as%
    \begin{equation*}
    g_\varepsilon(u) = 
\begin{cases}
\frac{1}{2}u^2 - u + 1, &  u \leq 1 + \varepsilon\\
\frac{1}{2}u^2 - \varepsilon \exp\left(\frac{u - 1}{\varepsilon} - 1\right), & 1 + \varepsilon < u < 1 + \varepsilon(1+\log 2)\\
\frac{1}{2}u^2 - 2u + 2 + 2\varepsilon \log 2, & u \geq 1 + \varepsilon(1+\log 2)
\end{cases}
\end{equation*}
for any $\varepsilon\geq 0$. Indeed, for $\varepsilon>0$, the objective is strongly convex and its unique critical point is $\bar x_{\varepsilon} = \exp\left(\frac{u-1}{\varepsilon}-1\right)$ so that the formula for $g_{\varepsilon}$ derives from verifying when $\bar x_{\varepsilon}\in[1,2]$. For $\varepsilon=0$ it suffices to check which vertex is optimal.

\smallskip
\ \ \ Observe that $u\in\mathbb R\mapsto \frac{1}{2}u^2-u+1$ achieves its minimum of $\frac 12$ at $\bar u_1=1\leq 1+\varepsilon$ whereas $u\in\mathbb R\mapsto \frac{1}{2}u^2-2u+2+2\varepsilon \log(2)$ achieves its minimum of $2\varepsilon\log(2)$ at $\bar u_2=2\geq 1+\varepsilon(1+\log(2))$ for $\varepsilon\leq \frac{1}{1+\log(2)}$. Conclude that, for each $0<\varepsilon\leq \frac{1}{1+\log(2)}$, $\bar u_1=1$ is a local minimizer of $g_{\varepsilon}$. However, $\bar u_1$ is merely a saddle point of $g_0$ as $g_0(u_-)>g_0(\bar u_1)>g_0(u_+)$ for every $u_-<\bar u_1<u_+<2$.%
\end{example}
\end{minipage}\hfill
\begin{minipage}{.35\textwidth}
      \begin{center} 
       \includegraphics[width=\linewidth]{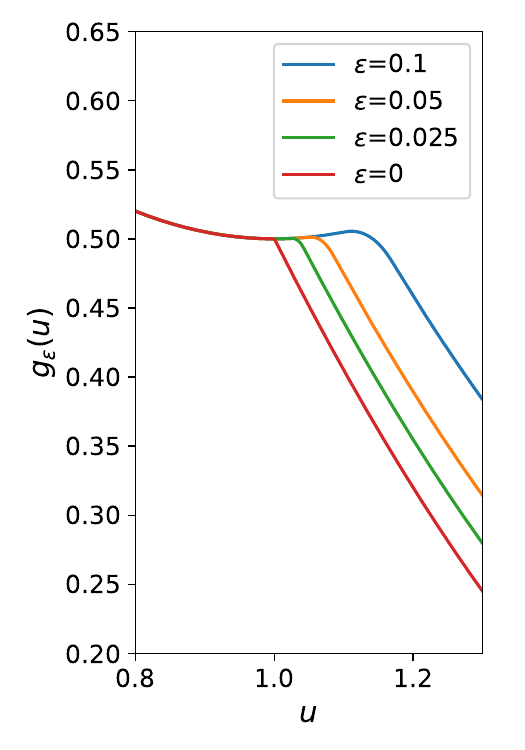}
       \vspace{-7.3mm}
       \captionsetup{width=0.9\linewidth}
          \noindent\captionof{figure}{Plot of $g_\varepsilon(u)$ for different values of $\varepsilon$. For any $\varepsilon>0$, $u=1$ is a local minimum of $g_\varepsilon$, but is a saddle~of~$g_0$. 
          \label{fig:1}}
       \end{center}
\end{minipage}

\vspace{4mm}

To address the issue identified in \cref{ex:convergencecritical}, we introduce \cref{alg:1} which, given a $\delta$-critical point $u_{\varepsilon}$ of $g_{\varepsilon}$, identifies a candidate for a local minimizer of $g_0$, $\bar u$, and certifies if $\bar u$ is a local minimizer of $g_0$ or not. 
We now describe the main routines in \cref{alg:1}.

\begin{figure}[!tb]
    \centering
    \includegraphics[width=0.7\linewidth]{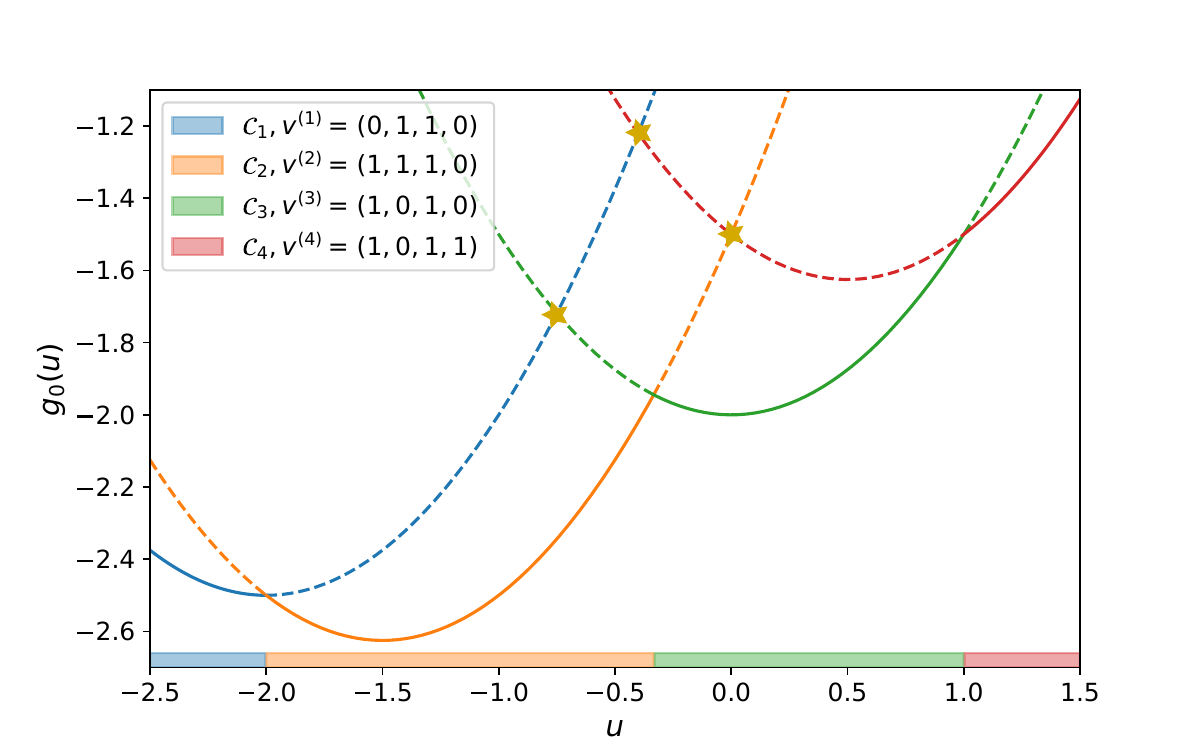}
    \caption{
    \label{fig:boundaryIllustration}
    Plot of the objective in the unregularized variational formulation for $\mathrm{M}=zz^{\intercal}$ with $z=(.5,-1.5,-.5,.5)$, $c=(-1,.5,-1,.5)$, and $\mathcal K=[0,1]^4$. Solid lines indicate the profile of $g_0$ and dashed lines are the plots of $g_0^{(i)}=\frac 12 u^2+(c-\rB^{\intercal}u)^{\intercal}v^{(i)}$ for $i\in[4]$. Stars indicate points where $g_0^{(i)}$ and $g_0^{(j)}$ meet when $\cC_i$ and $\cC_j$ do not share a common boundary. Such points lie above the graph of $g_0$.}
\end{figure}

\begin{itemize}
    \item \textbf{Candidate local minimizer of $\bm{g}_0$ (line 1):} First, we obtain the candidate point $\bar u=\rB \bar x$ for $\bar x \in \argmin_{x\in\mathcal K}(c-\rB^{\intercal}u_{\varepsilon})^{\intercal} x$. This choice is motivated by \cref{prop:nearCriticalPoints}, which asserts that there exist $\delta$-critical points of $g_{\varepsilon}$ that approximate local minimizers of $g_0$ well; namely, those $\delta$-critical points that are contained in the interior of some $\cC_i$ for which $i\in [V]$. 
    Thus, if such a $\cC_i$ contains a local minimizer of $g_0$, $\bar u$  is locally minimal by \cref{prop:minimizersG0}. 

\medskip
\item \textbf{Disqualify points violating necessary condition (lines 2-3):} Next, \cref{alg:1} discards $\bar u$ if  $(c-\rB^{\intercal}\bar u)^{\intercal} \bar x\neq\min_{x\in\mathcal K}(c-\rB^{\intercal}\bar u)^{\intercal} x$. This condition stems from \cref{prop:minimizersG0}, which implies that all local minimizers of $g_0$ satisfy $\bar u =\rB x_{\bar u}$ for every $x_{\bar u}\in\argmin_{x\in\mathcal K}(c-\rB^{\intercal}\bar u)^{\intercal} x$. Consequently, at this stage, the algorithm only discards points which are guaranteed to not be local minimizers. In particular, if $\bar u\in\Int(\cC_i)$ and $\cC_i$ does not contain $\rB v^{(i)}$, the algorithm will halt.  

\medskip
\item \textbf{Interior versus boundary point test (lines 4-18):} The remainder of \cref{alg:1} determines whether $\bar u$ is an interior or a boundary point of $\cC_i$. In the former case, $\bar u$ is a local minimizer of $g_0$, and is not otherwise. %
To start, we sample $w$ uniformly from the unit sphere, and fix the points $u_+=u+\eta w$ and $u_-=u-\eta w$ (lines 4-6), for some suitable $\eta>0$ that is obtained via a line search (lines 17-18). The algorithm proceeds as follows:

\medskip
\begin{itemize}[label=\raisebox{1.5pt}{\scalebox{.7}{$\blacktriangleright$}}]    

\item \textbf{Check if interior point (lines 9-11)}: It is shown in the proof of \cref{prop:alg_good} ahead, precisely in \cref{lem:interiorIFF}, that $\bar u \in \Int(\cC_i)$ if and only if,  for almost every choice of $w$ sampled uniformly from the unit sphere, there exists some $\eta>0$ for which $\rB x_+=\rB x_-$ for $x_+\in\argmin_{x\in\mathcal K}(c-\rB^{\intercal}u_+)^{\intercal}x$ and  $x_-\in\argmin_{x\in\mathcal K}(c-\rB^{\intercal}u_-)^{\intercal}x$. It is also shown that the values of $s_+,s_-,$ and $s$ do not coincide in this setting if $\rB x_+\neq\rB x_-$ so that the algorithm correctly identifies if $\bar u$ as an interior point and hence a local minimizer. %

\medskip  
\item \textbf{Check if boundary point (lines 13-16)}: In light of the above discussion, the algorithm halts correctly provided that the equality of $s_+,s_-,$ and $s$ implies that $\bar u$ is a boundary point as illustrated in 
\cref{fig:boundaryIllustration} and that this condition is eventually met. To this end it is shown in the proof of \cref{prop:alg_good} that if $\bar u\in\bd(\cC_i)$, $u_+$ and $u_-$ lie in distinct sets $\cC_+,\cC_-\in\{\cC_i\}_{i\in[V]}$ for every $\eta>0$ with probability $1$ (see \cref{lem:boundaryCondition1}). With this, once $\eta$ is sufficiently small that $\bar u\in\cC_+\cap \cC_-$, $u_+\in\Int(\cC_+)$, and $u_-\in\Int(\cC_-)$ the condition $s_+=s_-=s$ is shown to hold. \cref{lem:boundaryPointLemma} further guarantees that if $s_+=s_-=s$, $\bar u$ is a boundary point. 

\end{itemize}
\end{itemize}

\begin{algorithm}[!tb]
\caption{Certification}
\begin{algorithmic}[1]
\Statex \textbf{Input:} $u_\varepsilon \in \mathbb{R}^r$ and $\eta_0> 0$
\State $\bar u \gets \rB \bar x$ for $\bar{x} \in \argmin_{x \in \mathcal{K}} (c-\rB^\intercal u_\varepsilon )^\intercal x$ 
\If{$(c-\rB^\intercal \bar u )^\intercal \bar x > \min_{x \in \mathcal{K}} (c-\rB^\intercal \bar u )^\intercal x $}
    \State \Return $\bar u$ is not a local minimizer %
\EndIf
\State Sample $w$ uniformly from the unit sphere in $\mathbb R^r$
\State $\eta\gets \eta_0$
\State $u_{+}\gets \bar u +\eta w,\;u_{-}\gets \bar u -\eta w$ %
\State $s \gets \frac{1}{2} \|\bar u\|^2 - 
\min_{x \in \mathcal{K}}  (c-\rB^\intercal \bar u )^\intercal x$
\While{\textbf{True}}
    \State Pick $x_{+}\in \argmin_{x \in \mathcal{K}}  (c-\rB^\intercal u_{+})^\intercal x $, $x_{-}  \in \displaystyle{\argmin_{x \in \mathcal{K}}}  (c-\rB^\intercal u_{-} )^\intercal x$
    \If{$\rB x_{+}= \rB x_{-} $} %
        \State \Return $\bar u$ is  a local minimizer %
    \Else
        \State $s_{+}\gets \frac{1}{2} \|\bar u\|^2 - (\rB^\intercal \bar u - c)^\intercal x_{+}$, $s_{-}\gets \frac{1}{2} \|\bar u\|^2 - (\rB^\intercal \bar u - c)^\intercal x_{-} $
        \If{$s_{+}= s$\textbf{ and }$s_{-}= s$} 
            \State \Return $\bar u$ is not a local minimizer %
        \Else
           \State $\eta\gets \eta/2$ %
           \State $u_+\gets \bar u+\eta w,$ $u_-\gets \bar u-\eta w$ 
        
        \EndIf
    \EndIf
\EndWhile
\end{algorithmic}
\label{alg:1}
\end{algorithm}

\begin{figure}
    \centering
    \begin{subfigure}[b]{0.435\textwidth}
        \centering
        \includegraphics[width=0.97\textwidth]{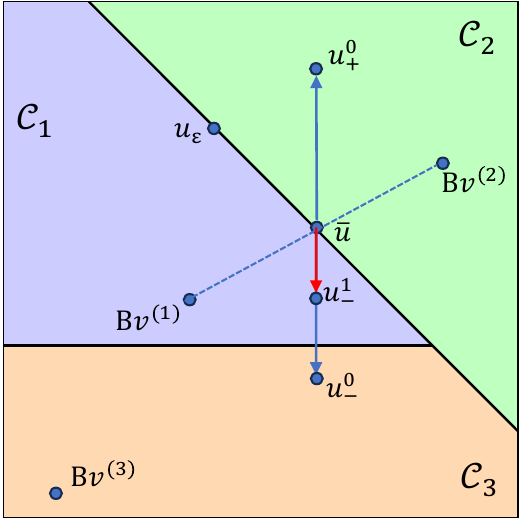}
        \caption{$\bar u \in \bd(\cC_1)$}
        \label{fig:bdry}
    \end{subfigure}
    \hfill
    \begin{subfigure}[b]{0.45\textwidth}
        \centering
        \includegraphics[width=\textwidth]{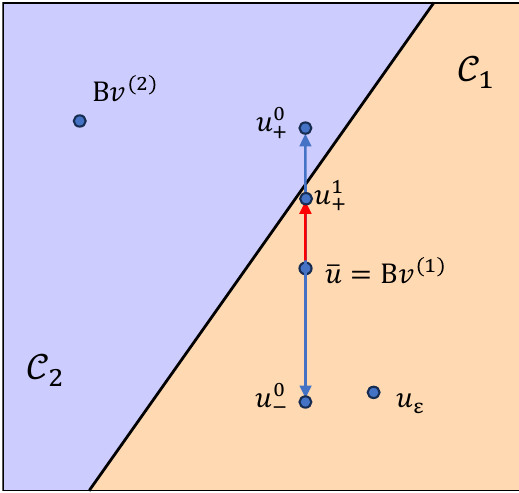}
        \caption{$\bar u \in \Int(\cC_1)$}
        \label{fig:interior}
    \end{subfigure}
    \caption{Visualization of \cref{alg:1} functioning when (a)  $\bar u$ lies on the boundary of $\cC_1$ and (b) $\bar u$ lies in the interior of $\cC_1$. The initially generated $u_+$ and $u_-$ are denoted by $u^0_+$ and $u^0_-$ respectively, while the perturbed vectors generated after $\eta$ is halved are denoted by $u^1_+$ and $u^1_-.$
    }
    \label{fig:combined}
\end{figure}

The following proposition establishes that the algorithm terminates correctly and provides an upper bound on the number of iterations of the while loop required. Its proof follows the main steps described above and is included in \cref{proof:prop:alg_good}. %

\begin{proposition}[Correctness and termination]\label{prop:alg_good} 
For $u\in \mathbb R^r$, let
$
\lambda(u) = \dist\left(u, \cup_{\{j \in [V] : u \notin \mathcal C_j\}} \mathcal C_j\right),$\text{ and } $
\zeta(u) = \frac{1}{2} \dist(u,\bd \cC_i).$  Then, \cref{alg:1} correctly 
identifies if $\bar u$ is a local minimizer or not  with probability 1 and 
\begin{enumerate}
        \item  if $\bar u \in \Int(\mathcal C_i)$,
         the while loop terminates within $\max\left\{0,\lceil \log_2\left({ \eta_0}/{\zeta(\bar u)}\right)  \rceil\right\}+1$ iterations.
        \item if $\bar u \in \bd(\mathcal C_i),$   the while  loop terminates within $\max\left\{0,\lceil \log_2\left({ \eta_0}/{\lambda(\bar u)}\right)  \rceil\right\}+2$ iterations.
    \end{enumerate}
    Furthermore, for any input $u_{\varepsilon}$, \cref{alg:1} terminates after at most 
    \[
    \max_{i \in [V]}\left\{
    \max\left\{0,\max_{u \in S \cap \Int (\mathcal C_i)}\lceil \log_2\left({ \eta_0}/{\zeta(u)}\right)  \rceil\right\}+1,
    \max\left\{0,\max_{u \in S \cap\bd (\mathcal C_i) }\lceil \log_2\left({ \eta_0}/{\lambda(u)}\right)  \rceil\right\}+2
     \right\} < \infty
    \] iterations, where $S$ is the set of Clarke critical points of $g_0$ which has finite cardinality.
\end{proposition}

\begin{remark}[Implementing \cref{alg:1}]
    Beyond simple arithmetic, \cref{alg:1} requires two main steps, namely, uniform sampling from the unit sphere in $\mathbb R^r$ and solving LPs.
    
To sample uniformly from the unit sphere, it suffices to first generate a sample from the standard Gaussian distribution in $\mathbb R^r$ and then to normalize this sample (see Problem 20.26 in \cite{billingsley2012probability}).   

Libraries for solving LPs, both commercial and open source, are plentiful. The vast majority are based on the simplex method, which has exponential worst-case complexity, or the interior point method, which has polynomial worst-case complexity. On most  practical examples, both approaches perform similarly.
\end{remark}

\section{Slow Rates Without Concavity}
\label{sec:slowRates}
While \cref{thm:EntropicQP} established an exponentially decaying entropic approximation gap for concave QPs, we next demonstrate that this fast rate is not attainable for convex or indefinite QPs. We start from a $O(\varepsilon)$ upper bound on the difference of (regularized and unregularized) costs, which follows from an elementary argument analogous to that used for LPs \cite{altschuler2017near,fang1993linear,weed2018explicit}, see \cref{proof:prop:slowRate} for details. %

\begin{proposition}
    \label{prop:slowRate}
    For each $\varepsilon>0$, let $x^{\star}_{\varepsilon}$ be a  solution of \eqref{eq:entropicQP} with regularization strength $\varepsilon$. Then, 
    \[
    \frac 12 \left(x^{\star}_{\varepsilon}\right)^{\intercal} \rM x^{\star}_{\varepsilon} + c^{\intercal} x^{\star}_{\varepsilon}- \min_{x\in\mathcal K}\left\{\frac 12 x^{\intercal} \rM x + c^{\intercal} x\right\}=O(\varepsilon).
    \]
    Moreover, if $x^{\star}_{\varepsilon}\to x^{\star}$ for some $x^{\star}\in \argmin_{x\in\mathcal K}\left\{\frac 12 x^{\intercal} \rM x + c^{\intercal} x\right\}$, this difference is of order $o(\varepsilon)$.
\end{proposition}

\subsection{The Convex Case}\label{sec:convexCounterexample}

We show that for convex QPs, no rate faster than $\varepsilon^2$ is attainable in general. To that end, consider the following example of a convex QP and its entropic regularization:
\[
    (\mathrm{P}_0)\coloneqq \min_{x\in[0,2]}\left\{\frac{1}{2}x^2-x\right\},\quad (\mathrm{P}_{\varepsilon})\coloneqq \min_{x\in[0,2]}\left\{\frac{1}{2}x^2-x+\varepsilon x\log x\right\}.
\]
It is easy to see that the global minimizer of $(\mathrm P_0)$ is given by $x^{\star}=1$ which achieves the value $-\frac 12$. On the other hand, the critical points of the objective in $(\mathrm P_{\varepsilon})$ consist of all points $\bar x_{\varepsilon}$ satisfying 
\[\bar x_{\varepsilon}e^{
    \frac{\bar x_{\varepsilon}}\varepsilon} =e^{\frac{1}{\varepsilon} -1},\text{ that is, } \bar x_{\varepsilon} = \varepsilon W_0\left( \varepsilon^{-1}{e^{\frac{1}{\varepsilon} -1}}\right), 
\]
where $W_0$ is the principal branch of the Lambert $W$ function.\footnote{For any $z\geq 0$, $W_0(z)$ is defined to be the unique solution, $y$, to the equation $ye^y=z$, see \cite{corless1996lambert} for details.} Equation (4.19) in \cite{corless1996lambert} yields 
\[
    \lim_{z\to \infty}
   \frac{\log(z)}{\log\left(\log(z)\right)}\left(W_0(z)-\log(z)+\log\left(\log(z)\right)\right)
      =1,
\]
and by substituting $z=\varepsilon^{-1}e^{\frac{1}{\varepsilon}-1}$, we further obtain
\[
    W_0\left(\varepsilon^{-1}{e^{\frac{1}{\varepsilon} -1}}\right)-\frac{1}{\varepsilon}+1-\log\left(\frac 1\varepsilon\right)+\log\left(\frac{1}{\varepsilon}-1+\log\left(\frac 1 \varepsilon\right)\right) = R(\varepsilon), 
\]
where $R(\varepsilon) = \Omega\left(\frac{\log\left(\frac{1}{\varepsilon}-1-\log(\varepsilon)\right)}{\frac{1}{\varepsilon}-1-\log(\varepsilon)}\right)$ or, equivalently, $ \Omega\left( \varepsilon\log\left(\frac 1\varepsilon\right)\right)$ as $\varepsilon\downarrow 0$. Notably,
\[
\begin{aligned}
   \frac{\varepsilon W_0\left(\varepsilon^{-1}{e^{\frac{1}{\varepsilon} -1} }\right)-{1}}{\varepsilon} &=-1-\log\left({1}-\varepsilon+\varepsilon\log\left(\frac 1\varepsilon\right)\right)+ R(\varepsilon)\to -1 \text{ as }\varepsilon\downarrow  0. 
\end{aligned}
\]
Thus, for every sufficiently small $\varepsilon$, $\varepsilon W_0\left(\varepsilon^{-1}{e^{\frac{1}{\varepsilon} -1} }\right)=\bar x_{\varepsilon}\in[0,2]$ and hence $\bar x_{\varepsilon}$ is the unique global minimizer of $(\mathrm{P}_{\varepsilon})$ over $[0,2]$. As $x^{\star}=1$,  ${\bar x_{\varepsilon}-x^{\star}}=\Omega(\varepsilon)$. The difference in costs is given as
\[
\begin{aligned}
\frac{1}{2}\bar x^2_{\varepsilon}-\bar x_{\varepsilon}-
    \min_{x\in[0,2]}\left\{\frac{1}{2}x^2-x\right\}&=  \frac{1}{2}(\bar x_{\varepsilon}-x^{\star})(\bar x_{\varepsilon}+x^{\star}) + x^{\star}-\bar x_{\varepsilon} = \frac{1}{2}(\bar x_{\varepsilon}-x^{\star})(\bar x_{\varepsilon}+x^{\star}-2),  
    \end{aligned}
\]
where this final term simplifies to $\frac{1}{2}(\bar x_{\varepsilon}-x^{\star})^2=\Omega(\varepsilon^2)$.

\subsection{The Indefinite Case} For indefinite QPs, i.e., when the objective is neither convex nor concave, we show that the entropic approximation gap is at best of order $\varepsilon^2$ in the absence of additional assumptions. %
Consider the problems 
\[
    (\mathrm{P}_0)\coloneqq \min_{\substack{x\in[0,2]\\y\in[0,2]}}\left\{\frac{1}{2}x^2-x-\frac{1}{2}y^2-y\right\},\quad (\mathrm{P}_{\varepsilon})\coloneqq     \min_{\substack{x\in[0,2]\\y\in[0,2]}}\left\{\frac{1}{2}x^2-x-\frac{1}{2}y^2-y+\varepsilon x\log(x)+\varepsilon y\log(y)\right\}.
\]
Clearly $(\mathrm P_0)$ is convex in $x$ and concave in $y$ so that the overall problem is neither convex nor concave. Furthermore, the problems are separable in the sense that 
\[
(\mathrm{P}_{\varepsilon}) = 
\min_{{x\in[0,2]\\}}\left\{\frac{1}{2}x^2-x+\varepsilon x\log(x)\right\}+\min_{y\in[0,2]}\left\{-\frac{1}{2}y^2-y +\varepsilon y\log(y)\right\},
\] 
where we highlight that the first problem coincides with the one studied in \cref{sec:convexCounterexample}. It is easy to see that the function $y\in[0,2]\mapsto -\frac{1}{2}y^2-y+\varepsilon y\log(y)$ is monotonically decreasing for every $0\leq \varepsilon\leq 1$ so that the minimum for the second problem is achieved at $\bar y_{\varepsilon}=2$. Conclude that if $(\bar x_{\varepsilon},\bar y_{\varepsilon})$ solves $(\mathrm{P}_{\varepsilon})$ and $(x^{\star}, y^{\star})$ solves $(\mathrm{P}_{0})$, it holds that
$
    \|(\bar x_{\varepsilon},\bar y_{\varepsilon})-(x^{\star}, y^{\star})\| = \left|\bar x_{\varepsilon}-x^{\star}\right| = \Omega(\varepsilon),
$
and that, for every $0\leq \varepsilon \leq 1$,  
\[
    \frac{1}{2}\bar x^2_{\varepsilon}-\bar x_{\varepsilon}-\frac{1}{2}\bar y^2_{\varepsilon}-\bar y_{\varepsilon}-\min_{\substack{x\in[0,2]\\y\in[0,2]}}\left\{\frac{1}{2}x^2-x-\frac{1}{2}y^2-y\right\}=\frac{1}{2}\bar x^2_{\varepsilon}-\bar x_{\varepsilon}-\min_{\substack{x\in[0,2]}}\left\{\frac{1}{2}x^2-x\right\}=\Omega(\varepsilon^2).
\]

\section{Duality and Approximation Rates for Gromov--Wasserstein Alignment}
\label{sec:GW}

The Gromov-Wasserstein problem provides a general framework for aligning metric measure space, with respect to a prescribed intrinsic notion of similarity on each space (e.g., the metric). Solving the GW problem not only quantifies how similar the metric measure spaces are, but also provides a scheme by which to optimally align them. The precise formulation is as follow; see  \cref{sec:primerOT} for background on the OT and GW problems.

Given finite sets $\mathcal X_0,\mathcal X_1$ with  probability measures $\mu_0,\mu_1$ on them    and some similarity measures $\kappa_{0}:\mathcal X_0\times \mathcal X_0\mapsto \mathbb R$ and  $\kappa_{1}:\mathcal X_1\times \mathcal X_1\mapsto \mathbb R$, the $(p,q)$-Gromov-Wasserstein problem with $p,q>0$ is  
\[
\mathsf{GW}_{p,q}(\mu_0,\mu_1)=\left(\inf_{\pi\in\Pi(\mu_0,\mu_1)} \iint \left|k_0^q(x,x')-k_1^q(y,y')\right|^p d\pi(x,y)\pi(x',y')\right)^{1/p}.
\]
Throughout, we assume for simplicity that $k_0$ and $k_1$ are symmetric functions. The EGW problem penalizes the alignment cost by the negative entropy and reads as
\begin{equation}
\label{eq:EGWFormula}
\mathsf{EGW}_{p,q}^\varepsilon(\mu_0,\mu_1)=\inf_{\pi\in\Pi(\mu_0,\mu_1)} \iint \left|k_0^q(x,x')-k_1^q(y,y')\right|^p d\pi(x,y) \pi(x',y')-\varepsilon\mathsf H(\pi),\quad \varepsilon>0,
\end{equation}
where $\mathsf H(\pi) =- \sum_{(x_0,x_1)\in\mathcal X_0\times \mathcal X_1} \pi\left(\left\{\left(x_0,x_1\right)\right\}\right)\log\left(\pi\left(\left\{\left(x_0,x_1\right)\right\}\right)\right)$. 

EGW was introduced in \cite{peyre2016gromov,solomon2016entropic} as a means to reduce the computational burden of solving the GW problem
directly. While the regularized problem may not be convex, \cite{solomon2016entropic} propose to optimize it using Sinkhorn iterations. This method is known to converge to a stationary point of a tight (albeit non-convex) relaxation of the EGW problem in the asymptotic regime, so that the overall computational complexity of this approach is unknown. Similar limitations apply for the popular mirror descent-based approach from \cite{peyre2016gromov} and its low-rank variant \cite{scetbon2022linear}. As noted in \cref{rem:entropic_computation}, \cite{rioux2023entropic} proposed an alternative method based on optimizing the variational form subject to non-asymptotic convergence guarantees, though this approach applies only  when $p=q=2$, $\mu_0,\mu_1$ are supported in Euclidean spaces, and $k_0,k_1$ are the corresponding Euclidean distances or the inner product.

We henceforth restrict our attention to the $p=2$ case, as the subsequent results only hold in this setting. When $\mu_0$ and $\mu_1$ are supported on a finite number of points, $(x_0^{(i)})_{i=1}^{N_0}\subset\mathcal X_0$ and $(x_1^{(i)})_{i=1}^{N_1}\subset\mathcal X_1$, respectively, the entropic GW problem can be recast as a QP:
\begin{equation}
\label{eq:finiteEGWDecomposition}
\begin{aligned}
\mathsf{EGW}_{2,q}^\varepsilon(\mu_0,\mu_1)= \iint k_0^{2q}d\mu_0d \mu_0+ \iint k_1^{2q}d\mu_1d\mu_1+ \min_{\substack{\mathrm A x = b\\x\geq 0}}\left\{ -\frac 12 x^{\intercal}\mathrm{M}x - \varepsilon \mathsf H(x)\right\},
\end{aligned}
\end{equation}
where 
 the constraint $\pi\in\Pi(\mu_0,\mu_1)$ is rewritten as $x\in\{x\in\mathbb R^{N_0N_1}:\mathrm {A}x=b,x\geq 0\}$ for  
\begin{equation}
\label{eq:Ab}
\mathrm{A} = \begin{pmatrix}
   \mathbf{1}_{N_1}^{\intercal}&\mathbf 0_{N_1}^{\intercal}&\cdots&\mathbf 0_{N_1}^{\intercal}\\ 
   \mathbf 0_{N_1}^{\intercal}&\mathbf{1}_{N_1}^{\intercal}&\cdots&\mathbf 0_{N_1}^{\intercal}\\
   \vdots&\vdots &\ddots&\vdots\\
   \mathbf 0_{N_1}^{\intercal}&\mathbf 0_{N_1}^{\intercal}&\cdots&\mathbf{1}_{N_1}^{\intercal}&\\
   \mathrm{Id}_{N_1}&\mathrm{Id}_{N_1}&\cdots&\mathrm{Id}_{N_1}
\end{pmatrix}\in \mathbb R^{(N_0+N_1)\times N_0N_1},\quad b = \begin{pmatrix}
   \mu_0\left(\left\{x^{(1)}_0\right\}\right)\\\vdots\\ \mu_0\left(\left\{x^{(N_0)}_0\right\}\right)\\ \mu_1\left(\left\{x^{(1)}_1\right\}\right)\\\vdots\\ \mu_1\left(\left\{x^{(N_1)}_1\right\}\right)
\end{pmatrix}\in\mathbb R^{N_0+N_1},
\end{equation}
with $\mathbf{1}_{N_1}$ and $\mathbf{0}_{N_1}$ as the vectors of all $1$'s and all $0$'s of length $N_1$, $\mathrm{Id}_{N_1}\in\mathbb R^{N_1\times N_1}$ the identity matrix, and 
\begin{equation}
\label{eq:M}
\mathrm{M} = 4\left(\kappa _0^q\left(x^{(i)}_0,x^{(j)}_0\right)\right)_{i,j=1}^{N_0}\otimes \left(\kappa _1^q\left(x^{(i)}_1,x^{(j)}_1\right)\right)_{i,j=1}^{N_1}\in\mathbb R^{N_0N_1\times N_0N_1}, 
\end{equation}
where $\otimes$ denotes the Kronecker product of two matrices.  We highlight that the constraints $\pi\in\Pi(\mu_0,\mu_1)$ and $x\in\{x\in\mathbb R^{N_0N_1}:\mathrm {A}x=b,x\geq 0\}$ can be seen to be equivalent under the transformation $x= \vec(\rP)$  where $\rP\in\mathbb R^{N_1\times N_0}$ has entries $\rP_{lm}=\pi\left(\left\{\left(x_0^{(m)},x_1^{(l)}\right)\right\}\right)$ and $\vec(\rP)$ is the vector obtained by stacking the columns of $\rP$.

 With \eqref{eq:finiteEGWDecomposition} in hand, the following result is an immediate consequence of \cref{prop:variationalFormND} and the properties of PSD and negative semidefinite (NSD) kernels.\footnote{ A PSD (resp. NSD) kernel $k$ is a symmetric function for which, for any  $s\in\mathbb N$ and $\bm x\coloneqq (x_1,\dots,x_s)\in\prod_{i=1}^s\mathcal X_i$, the matrix $\mathrm E_{\bm x}\in\mathbb R^{s\times s}$ with $lm$-th entry given by $(\mathrm E_{\bm x})_{lm}=k(x_l,x_m)$ is PSD (resp. NSD).} 

\begin{proposition}[EGW variational form]
\label{prop:EGWDuality}
Let $(\mathcal X_0,\mu_0)$, $(\mathcal X_1,\mu_1)$ be finite probability spaces and $\kappa_{0}:\mathcal X_0\times \mathcal X_0\to \mathbb R$, $\kappa_{1}:\mathcal X_1\times \mathcal X_1\to \mathbb R$ be both PSD or both NSD kernels. Setting $\mathsf C({\mu_0,\mu_1})=\iint \kappa_{0}^{2  q}d\mu_0d\mu_0+\iint \kappa_{1}^{2q}d\mu_1d \mu_1$, for $q\in\mathbb N\backslash \{0\}$, we have that \[
\mathsf{EGW}^{\varepsilon}_{2,q}(\mu_0,\mu_1)
  =\mathsf C({\mu_0,\mu_1})+ \inf_{u\in\mathbb R^{r}}\left\{\frac{1}{2}\|u\|^2+\min_{\substack{\mathrm{A}x=b\\x\geq 0}}\left\{- u^{\intercal}\rB x-\varepsilon \mathsf H(x)  \right\}\right\},
\]
for each $\varepsilon\geq 0$,
where $\rB \in\mathbb R^{r\times N_0N_1}$  satisfies $\mathrm {M}=\mathrm{B}^{\intercal}\mathrm{B}$ and $r\leq N_0N_1$. 

\end{proposition}

The claim follows from \eqref{eq:finiteEGWDecomposition} by applying  \cref{prop:variationalFormND}. To invoke the proposition, we note that the positive/negative semidefiniteness of the two kernels ensures that the matrix $\mathrm{M}$ is PSD. Indeed, $\mathrm{M}$ is then the Kronecker product of PSD/NSD matrices  and is thus PSD by Schur's product theorem \cite{schur1911}.\footnote{If $\mathrm{R}$ is a PSD matrix, the Hadamard product $\mathrm R\odot \mathrm R$ (which squares all entries of $\mathrm{R}$) is also PSD; cf. e.g., Theorem 4.2.12 in \cite{horn1991topics}. If $\mathrm{S}$ is NSD, the same result implies that $\mathrm S\odot \mathrm{S}$ is PSD. It follows that the matrices in \eqref{eq:M} are both PSD/NSD, depending on $q$.}

\begin{remark}[Comparison to variational form in \cite{zhang2024gromov}]\label{rem:2-2}

   The recent work \cite{zhang2024gromov} derived a similar variational form for $(2,2)$-GW distances between  centered distributions on Euclidean spaces with finite fourth moments with and without regularization (centering does not restrict generality as GW and EGW are invariant to translations).  This corresponds to the case where $\mathcal X_i=\mathbb R^{d_i}$ and $\kappa_i(x,x')=\|x-x'\|$ for $i\in\{0,1\}$. They demonstrate that $\mathsf{GW}_{2,2}^\varepsilon(\mu_0,\mu_1)=\mathsf K({\mu_0,\mu_1})+\mathsf R^{\varepsilon}(\mu_0,\mu_1)$, where 
\begin{equation}
\label{eq:quadEGW}
\begin{aligned}
&\mathsf K({\mu_0,\mu_1})\mspace{-4mu}=\mspace{-4mu}\iint\kappa_0^4d\mu_0 d\mu_0+\iint\kappa_1^4d\mu_1 d\mu_1-4\iint\|x\|^2\|y\|^2d\mu_0(x)d\mu_1(y),\\
&\mathsf R^{\varepsilon}(\mu_0,\mu_1)\mspace{-4mu}=\mspace{-4mu}\inf_{\pi\in\Pi(\mu_0,\mu_1)}\int\mspace{-3mu}-4\|x\|^2\|y\|^2d\pi(x,y)\mspace{-3mu}-\mspace{-3mu}8\mspace{-3.5mu}\sum_{\substack{1\leq i\leq d_0\\1\leq j\leq d_1}}\mspace{-4mu}\left(\int\mspace{-3mu} x_iy_j\,d\pi(x,y)\right)^2\mspace{-3mu}-\varepsilon\mathsf H(\pi),
\end{aligned}
\end{equation}
   for any $\varepsilon\geq 0$ and that $\mathsf{R}^{\varepsilon}$ admits the variational form 
   \begin{equation}
   \label{eq:varFormR}
       \mathsf{R}^{\varepsilon}(\mu_0,\mu_1) = \inf_{\mathrm Z\in\mathbb R^{d_0\times d_1}}\left\{32\|\mathrm Z\|_{\mathrm F}^2+\inf_{\pi\in\Pi(\mu_0,\mu_1)}\int -4\|x\|^2\|y\|^2-32x^{\intercal}\mathrm Zyd\pi(x,y)-\varepsilon\mathsf H(\pi)\right\},   
     \end{equation}
where $\|\cdot\|_{\mathrm{F}}$ denotes the Frobenius matrix norm.

On the other hand, it is not immediately clear that the conditions of \cref{prop:EGWDuality} hold when $\mu_0$ and $\mu_1$ are centered and finitely supported, as $\mathrm{M}$ in \eqref{eq:M} is not generally a positive semidefinite matrix for this choice of $k_0,k_1$. However, using the decomposition in \eqref{eq:quadEGW} and the correspondence between couplings and $\{x\in\mathbb R^{N_0N_1}:\mathrm{A}x=b,x\geq 0\}$, we 
see that 
\begin{equation}\label{eq:quad2-2}
    -\frac 12 
x^{\intercal} \mathrm{M}x = -\frac 12  x^{\intercal} \mathrm{M}'x+c^{\intercal} x+C
\end{equation}
where $C = -4\int \|x\|^2\|y\|^2d\mu_0(x) \mu_1(y)$, 
\[
c =-4\left(\left\|x_0^{(1)}\right\|^2\left\|x_1^{(1)}\right\|^2,\dots,\left\|x_0^{(1)}\right\|^2\left\|x_1^{(N_y)}\right\|^2,\dots,\left\|x_0^{(N_x)}\right\|^2\left\|x_1^{(1)}\right\|^2,\dots,\left\|x_0^{(N_x)}\right\|^2\left\|x_1^{(N_y)}\right\|^2 \right)^{\intercal},
\] 
and $\rM'=\rB^{\intercal}\rB$ where $
\rB^{\intercal}=4\left(\begin{smallmatrix}
   |&&|&&|&&| 
    \\z_{11}&\cdots&z_{1d_1}&\cdots&z_{d_01}&\cdots&z_{d_0d_1}
    \\
    |&&|&&|&&|
\end{smallmatrix}\right)$  for
\[
\begin{aligned} z_{ij}&=\left(\left(x_0^{(1)}\right)_i\left(x_1^{(1)}\right)_j,\dots,\left(x_0^{(1)}\right)_i\left(x_1^{(N_y)}\right)_j,\dots,\left(x_0^{(N_x)}\right)_i\left(x_1^{(1)}\right)_j,\dots,\left(x_0^{(1)}\right)_i\left(x_1^{(N_y)}\right)_j\right)^{\intercal}.
\end{aligned}
\] 
\cref{prop:variationalFormND} then yields the alternative, but equivalent, variational form for $\mathsf R^{\varepsilon}(\mu_0,\mu_1)$,
\[
\mathsf{R}^{\varepsilon}(\mu_0,\mu_1)=
\inf_{u\in\mathbb R^{d_0d_1}}\left\{\frac 12 \|u\|^2+\min_{\substack{\mathrm{A}x=b\\x\geq 0}}\left\{(c-\mathrm{B}^{\intercal}u)^{\intercal}x-\varepsilon\mathsf{H}(x)\right\}\right\}. 
\] 
To see the equivalence to \eqref{eq:varFormR}, it suffices to set $u=8\mathrm{vec}(Z)$ in the above problem.

A similar result for the inner product cost $\kappa_i:(x,x')\in\mathbb R^{d_i}\times \mathbb R^{d_i}\mapsto x^{\intercal}x'$ appeared in Appendix~E of \cite{rioux2023entropic} and can be recovered using the same logic as above.
\end{remark}

Beyond the two cases noted in \cref{rem:2-2}, no variational formulations for GW problems have appeared previously in the literature. As such, \cref{prop:EGWDuality} provides one for a large number of new settings. %
Noting that the variational form for the Euclidean quadratic GW distance enabled new analytic and computational techniques to be developed, we anticipate that the results above with unlock similar progress for discrete GW over general spaces and arbitrary PSD/NSD kernels.

\medskip
\cref{rem:2-2} enables leveraging \cref{thm:EntropicQP} to derive, for the first time, an exponential entropic approximation rate for the discrete GW problem, with a wide array of costs. %

\begin{corollary}
\label{cor:costGW}
     Let $\mathcal X_i$, $\kappa_i$, and $\mu_i$ be as in \cref{prop:EGWDuality} with $q> 0$ or \cref{rem:2-2} with $q=2$. Then, letting $x_{\varepsilon}^{\star}$ solve the regularized QP in \eqref{eq:finiteEGWDecomposition} for each $\varepsilon> 0$,   
     \[
        \dist\left(x^{\star}_{\varepsilon},\argmin_{\substack{\mathrm{A}x=b\\x\geq 0}}-\frac{1}{2}x^{\intercal}\rM x\right)\leq C_1\exp(-C_2/\varepsilon), 
     \]
     for every $0<\varepsilon\leq C_3$ where $C_1,C_2,$ and $C_3$ are given explicitly in \cref{thm:EntropicQP}. Moreover, if $\pi^{\star}_{\varepsilon}$ solves \eqref{eq:EGWFormula} for every $\varepsilon>0$,
     \[
       0\leq \iint \left|k_0^q(x,x')-k_1^q(y,y')\right|^2 d\pi_{\varepsilon}^{\star}(x,y)d \pi_{\varepsilon}^{\star}(x',y')-\mathsf{GW}_{2,q}^2(\mu_0,\mu_1)\leq C_1'\exp(-C_2/\varepsilon)  
     \]
     for every $0<\varepsilon\leq C_3$, where $C_2$ and $C_3$ are as above and $C_1'$ is as in \cref{cor:cost}.
\end{corollary}

\section{Proofs of Main Results}
\label{sec:proofs}
\subsection{Proofs for \texorpdfstring{\cref{sec:concaveQP}}{Section 3}}
\subsubsection{Proof of \texorpdfstring{\cref{prop:variationalFormND}}{Proposition 1} }
\label{proof:prop:variationalFormND}
 As $\rM$ is positive semidefinite, the function $f(x)=\frac{1}{2}x^{\intercal}\rM x-c^{\intercal} x$ is convex with convex conjugate
    \[
    f^{*}:y \in \mathbb R^n\mapsto \frac{1}{2}(y+c)^{\intercal}\rM^{\dagger}(y+c) + \delta_{\range(\rM)}(y+c)
    \]
    (see Example 11.10 in \cite{rockafellar2009variational}). It  
      follows from the Fenchel-Moreau theorem (cf. e.g., Theorem 4.2.1 in \cite{borwein2006convex}) that 
        \begin{equation}
        \label{eq:fenchelMoreau}
      \begin{aligned}
          f(x) &= \sup_{y\in\mathbb R^n}\left\{x^{\intercal}y-\frac{1}{2}(y+c)^{\intercal}\rM^{\dagger}(y+c)-\delta_{\range(\rM)}(y+c)\right\}\\
          &= - \inf_{y\in\mathbb R^n}\left\{\frac{1}{2}(y+c)^{\intercal}\rM^{\dagger}(y+c) - x^{\intercal}y +  \delta_{\range(\rM)}(y+c)\right\}.
      \end{aligned}
      \end{equation}
    so that 
    \[
        \min_{x\in\mathcal K}\left\{-f(x)-\varepsilon\mathsf H(x)\right\} = \inf_{y\in\mathbb R^n}\left\{\frac{1}{2}(y+c)^{\intercal}\rM^{\dagger}(y+c)+\inf_{x\in\mathcal K}\left\{-x^{\intercal}y-\varepsilon\mathsf H(x)\right\}+  \delta_{\range(\rM)}(y+c)\right\}.
    \]
    As $y+c\in\range(\rM)$, we may instead pass to an auxiliary variable $z\in\mathbb R^n$ and substitute $y=\rM z-c$ in the above expression, yielding 
    \[
    \begin{aligned}
\min_{x\in\mathcal K}\left\{-f(x)-\varepsilon\mathsf H(x)\right\} &= \inf_{z\in\mathbb R^n}\left\{\frac{1}{2}z^{\intercal}\rM z+\inf_{x\in\mathcal K}\left\{-x^{\intercal}(\rM z-c)-\varepsilon\mathsf H(x)\right\}\right\}
\\
&=
\inf_{u\in\rB\mathbb R^n}\left\{\frac{1}{2}\|u\|^2+\inf_{x\in\mathcal K}\left\{-x^{\intercal}(\rB^{\intercal} u-c)-\varepsilon\mathsf H(x)\right\}\right\}
    \end{aligned} 
    \]
    where, in the first equality, we have used the fact that $\rM\rM^{\dagger}\rM =\rM$ and, in the second, we have set $u=\rB z$. 

    We now show that all solutions to \eqref{eq:variationalFormND} are necessarily elements of $\rB\mathbb R^n$ so that the second infimum in the previous display can be taken over $\mathbb R^r$. First, note that, for any $x'\in\mathcal K$, 
    \begin{equation}
    \label{eq:gCoercive}
        g_{\varepsilon}(u) \coloneqq \frac{1}{2}\|u\|^2+\inf_{x\in\mathcal K}\left\{-x^{\intercal}(\rB^{\intercal} u-c)-\varepsilon\mathsf H(x)\right\}\geq \frac{1}{2}\|u\|^2 -\|\rB x'\|\|u\|+ c^{\intercal}x'-\varepsilon \mathsf{H}(x'),  
    \end{equation}
    so that $g_{\varepsilon}$ is coercive and hence each of its global minimizers over $\mathbb R^r$ is necessarily a local minimizer. It is easy to see that $g_{\varepsilon}$ is locally Lipschitz continuous\footnote{For any compact set $D\subset \mathbb R^r$ and $u,u'\in D$, $\left|\|u\|^2-\|u
'\|^2\right| \leq (\|u\|+\|u'\|)\|u-u'\|\leq 2\sup_{\tilde u\in D}\|\tilde u\|\|u-u'\|$ and $\left|\inf_{x\in\mathcal K}\left\{-x^{\intercal}(\rB^{\intercal}u-c)-\varepsilon\mathsf H(x)\right\}-\inf_{x\in\mathcal K}\left\{-x^{\intercal}(\rB^{\intercal}u'-c)-\varepsilon\mathsf H(x)\right\}\right|\leq \sup_{x\in\mathcal K}\left|(u'-u)^{\intercal}\rB x\right|\leq\sup_{x\in\mathcal K}\|\rB x\|\|u-u'\|$.} on $\mathbb R^r$ so that Proposition 2.3.2 in \cite{clarke1990optimization} asserts that a local minimizer, $\bar u$, of $g_{\varepsilon}$ must satisfy 
    \[
       0\in\partial g_{\varepsilon}(\bar u) = \bar u + \partial V(\bar u),\text{ where }V(u)\coloneqq \min_{x\in\mathcal K}\left\{-x^{\intercal}(\rB^{\intercal} u-c)-\varepsilon\mathsf H(x)\right\},  
    \]
    where we have applied the subdifferential sum rule (see Corollary on p. 39 of \cite{clarke1990optimization}). It is well known (cf. e.g., Proposition B.22 in \cite{bertsekas1999Nonlinear}) that $-V$ is convex and has subdifferential\footnote{While the cited result pertains to the subdifferential in the sense of convex analysis, this set coincides with the Clarke subdifferential by virtue of Proposition 2.2.7 in \cite{clarke1990optimization}.} given by 
    \[
    \begin{aligned}
        \partial(-V)(u) &= \conv\left(\left\{ \rB x:x\in\argmax_{x\in\mathcal K}\left\{x^{\intercal}(\rB^{\intercal} u-c)+\varepsilon\mathsf H(x)\right\}\right\}\right)
        \\
        &=  \rB \argmin_{x\in\mathcal K}\left\{-x^{\intercal}(\rB^{\intercal} u-c)-\varepsilon\mathsf H(x)\right\},
    \end{aligned} 
    \]
    where we have used the fact that the solution set of a convex problem is  convex. As $\partial(-V)=-\partial V$ by Proposition 2.3.1 in \cite{clarke1990optimization}, we must have that a global minimizer, $u^{\star}$, of $g_{\varepsilon}$ satisfies 
    \[
        u^{\star}\in \rB \argmin_{x\in\mathcal K}\left\{-x^{\intercal}(\rB^{\intercal} u^{\star}-c)-\varepsilon\mathsf H(x)\right\} \subset \rB \mathcal K\subset \rB \mathbb R^n, 
    \]
    proving the first claim. 

    As for points (1) and (2), if $x^{\star}\in\mathcal K$ is optimal for \eqref{eq:primalQPND}, then 
    \[
        \inf_{u\in\mathbb R^r} g_{\varepsilon}(u)\leq g_{\varepsilon}(\rB x^{\star}) \leq -\frac{1}{2} (x^{\star})^{\intercal}\rM x^{\star}+c^{\intercal}x^{\star}-\varepsilon \mathsf{H}(x^{\star}) = \inf_{u\in\mathbb R^r} g_{\varepsilon}(u),
    \]
   so that both inequalities are in fact equalities, i.e., both statements in point (1) hold. On the other hand,     
    for any $u\in\mathbb R^r$, let $x_u\in\argmin_{x\in\mathcal K}\left\{-x^{\intercal}(\rB^{\intercal} u-c)-\varepsilon\mathsf H(x)\right\}$ be arbitrary. Then, 
    \[
        g_{\varepsilon}(u) =\frac{1}{2}\|u-\rB x_u\|^2-\frac{1}{2}x^{\intercal}_u\rB^{\intercal}\rB x_{u} +x_u^{\intercal}c-\varepsilon \mathsf H(x_u)=\frac{1}{2}\|u-\rB x_u\|^2+f(x_{u})-\varepsilon \mathsf H(x_u), 
    \]
    so that 
   \[
        g_{\varepsilon}(u)-\inf_{u\in\mathbb R^r} g_{\varepsilon}(u) = \underbrace{f(x_u) - \varepsilon \mathsf H(x_u) -\inf_{x\in\mathcal K} \left\{f(x) - \varepsilon \mathsf H(x)\right\}}_{\geq 0} +\underbrace{\frac{1}{2}\|u-\rB x_{u}\|^2}_{\geq 0}.
   \]
   Thus, if $u^{\star}$ is optimal for \eqref{eq:variationalFormND}, we must have that both inequalities above are saturated, i.e., every element  $x_{u^{\star}}\in\argmin_{x\in\mathcal K}\left\{-x^{\intercal}(\rB^{\intercal} u^{\star}-c)-\varepsilon\mathsf H(x)\right\}$ solves \eqref{eq:primalQPND} and $u^{\star}=\rB x_{u^{\star}}$.
\hfill\qed

\subsubsection{Proof of \texorpdfstring{\cref{prop:minimizersG0}}{Proposition 2}}
\label{proof:prop:minimizersG0}

   Suppose that $u\in \mathcal C_i\cap \mathcal C_j$ for some distinct indices $i,j\in [V]$ so that $(c-\rB^{\intercal}u)^{\intercal}v^{(i)}=(c-\rB^{\intercal}u)^{\intercal}v^{(j)}$. If $\mathrm{B}^{\intercal} v^{(i)}=\mathrm{B}^{\intercal} v^{(j)}$, it follows that $c^{\intercal} v^{(i)}=c^{\intercal} v^{(j)}$ whereby $\mathcal C_i=\mathcal C_j$. Now, if $\mathrm{B} v^{(i)}\neq \mathrm{B} v^{(j)}$ and $u\in \Int(\mathcal C_i)$ (the argument follows similar lines with $\mathcal C_j$ in place of $\cC_i$), there exists $t>0$ sufficiently small that $u'=u-t\frac{\mathrm{B}v^{(i)}-\mathrm{B}v^{(j)}}{\|\mathrm{B}v^{(i)}-\mathrm{B}v^{(j)}\|}\in \mathrm{int}(\mathcal C_i)$, however, 
    \[
        (c-\rB^{\intercal} u')^{\intercal}(v^{(i)}-v^{(j)}) =  (c-\rB^{\intercal}u)^{\intercal}(v^{(j)}-v^{(i)})+ t\|\rB v^{(j)}-\rB v^{(i)}\| >0    \]
    so that $(c-\rB^{\intercal} u')^{\intercal}v^{(i)}>(c-\rB^{\intercal} u')^{\intercal}v^{(j)}$, which contradicts the fact that $u'\in \Int(\mathcal C_i)$. Conclude that, in this case, $u\in \bd(\mathcal C_i)$ and similarly for $\mathcal C_j$, proving the first assertion.

We now move to proving the statement on minimizers of $g_0$. Fix some $u\in \bd(\mathcal C_i)$ and assume,
without loss of generality,  that $i=1$ and that $u\in\cap_{j=1}^l\mathcal C_j$ and that this collection is maximal (i.e., $u\not\in C_k$ for any $k\in[V]\backslash[l]$). By definition of $\mathcal C_i$, 
    \begin{align}\label{step1}
        (c - \rB^\intercal u)^\intercal v^{(1)} = (c - \rB^\intercal u)^\intercal v^{(j)}\text{ for every $j\in[l]$ }.
    \end{align}
    Now, fix any $w \in \RR^r$ for which $u + w \in \Int(\cC_1)$ and note that, for every $j\in [l]\backslash \{1\}$,
    \begin{align}\label{step2}
        \frac{1}{2}\|u + w\|^2 + (c - \rB^\intercal (u + w))^\intercal v^{(1)} < \frac{1}{2}\|u + w\|^2 + (c - \rB^\intercal (u + w))^\intercal v^{(j)},
    \end{align}
    by definition of $\cC_1$. Together, 
    \eqref{step1} and \eqref{step2} yield that $w^\intercal (-\rB v^{(1)}) < w^\intercal (-\rB v^{(j)})$ whereby 
     \begin{equation}
     \label{eq:boundOnDirectionalDerivatives}\left(u - \rB v^{(1)}\right)^\intercal(-w) > \left(u - \rB v^{(j)}\right)^\intercal(-w) \text{ for every $j\in[l]\backslash\{1\}$}.\end{equation}
    As shown in   the proof of \cref{prop:variationalFormND}, $\partial g_0(u) = \conv\left(\{u - \rB v^{(i)} : i \in [l]\}\right)$ so that 
    \begin{align*}
        \limsup_{u \to u_0,t \downarrow 0}\frac{-g_0(u + tw) + g_0(u)}{t}= \max_{r \in \partial g_0(u)} r^\intercal(-w) = 
        \left(u - \rB v^{(1)}\right)^\intercal (-w)
    \end{align*}
    where the first equality follows from Propositions 2.1.1(c) and 2.1.2(b) in \cite{clarke1990optimization} and the second is due to \eqref{eq:boundOnDirectionalDerivatives}. We now show that the limit in the above display (called the generalized directional derivative) coincides with the the directional derivative of $g_0$ at $u$ in the direction $-w$. To this end, note that
    \begin{align*}
        -g_0(u) &= \max_{i \in [n]}\left\{-\frac{1}{2}\|u\|^2 + (\rB^\intercal - c)^\intercal v^{(i)}\right\},
    \end{align*}
    so that $-g_0$ is regular (in particular the generalized directional derivative and the standard directional derivative coincide, see Definition 2.3.4 in \cite{clarke1990optimization}) as follows by Proposition 2.3.12 in \cite{clarke1990optimization}. Conclude that 
    \begin{align*}
        \left(u - \rB v^{(1)}\right)^\intercal (-w) = \lim_{t \to 0}\frac{-g_0(u + tw) + g_0(u)}{t}.
    \end{align*}

    It remains to consider the sign of the directional derivative. 
    First, if it is positive for some $w\in \Int(C_1)-u$, we have that $\frac{-g_0(u + tw) + g_0(u)}{t}>0$ for every $t>0$ sufficiently small whereby $g_0(u+tw)<g_0(u)$ for all such $t$ as desired. Otherwise, the directional derivative is nonpositive for every $w\in\Int(C_1)-u$. %
    Thus, for any  $w\in\Int(C_1)-u$, we have that 
    \[
    \left(u-\rB v^{(j)}\right)^{\intercal}(-w)<\left(u-\rB v^{(1)}\right)^{\intercal}(-w)\leq  0
    \text{ for every $j\in[l]\backslash \{1\}$}
    \]
    as follows from \eqref{eq:boundOnDirectionalDerivatives}.  Applying the formulas from the first part of the proof,
    \[
        \lim_{t\downarrow 0} \frac{-g_0(u-tw)+g_0(u)}{t} = \max_{r\in\partial g_0(u)} r^{\intercal}w = \max_{\substack{\lambda_1,\dots,\lambda_l\geq 0\\\sum_{j=1}^l\lambda_j=1}}\sum_{j=1}^l\lambda_j\left(u-\rB v^{(j)}\right)^{\intercal}w>0,
    \]
    where the final equality is because $\partial g_0(u) = \conv\left(\{u - \rB v^{(i)} : i \in [l]\}\right)$ whereas the  inequality follows from the fact that   $\left(u-\rB v^{(j)}\right)^{\intercal}(-w) <  0$
    \text{ for every $j\in[l]\backslash \{1\}$}. Conclude that 
     $g_0(u)>g_0(u-tw)$ for every $t>0$ sufficiently small, proving the claim.

     This implies that every local minimizer, $\bar u$, of $g_0$ must be an interior point of  $\mathcal C_i$ for some $i\in[V]$. As $g_0(u)=\frac 12 \|u\|^2+(c-\rB^{\intercal}u)^{\intercal}v^{(i)}$ for every $u\in \mathcal C_i$, we see that $g_0$ attains a local minimum in $\Int(\mathcal C_i)$ if and only if $\rB^{\intercal}v^{(i)}\in \Int(\mathcal C_i)$ and this point corresponds to a local minimizer.  
\hfill\qed

\subsubsection{Proof of \texorpdfstring{\cref{thm:EntropicQP}}{Theorem 3}}
\label{proof:thm:entropicQP}
Fix $0 < \varepsilon < \min\left\{
  \frac{\alpha}{R_{\mathsf H}},
  \frac{\gamma^2}{8 R_{\mathsf H}},
  \frac{\Delta}{R_1 + R_{\sH}}
\right\}$ and let  $x_{\varepsilon}^{\star}$ solve the regularized QP \eqref{eq:primalQPND}. We first establish that $\rB x_{\varepsilon}^{\star} \in \mathcal C_i$ for some $i\in \cI^{\star} \subseteq [V]$.

To this end, fix $u\in\mathbb R^r$,  $x_{\varepsilon,u}\in\argmin_{x\in\mathcal K}\left\{(c-\rB^{\intercal}u)^{\intercal }x-\varepsilon\mathsf{H}(x)\right\}$, and set $x_{0,u}$ similarly for $\varepsilon=0$. Then, $ (c-\rB^{\intercal}u)^{\intercal }x_{0,u}\leq (c-\rB^{\intercal}u)^{\intercal }x_{\varepsilon,u}$ and $(c-\rB^{\intercal}u)^{\intercal }x_{\varepsilon,u}-\varepsilon\mathsf{H}(x_{\varepsilon,u})\leq (c-\rB^{\intercal}u)^{\intercal }x_{0,u}-\varepsilon\mathsf{H}(x_{0,u})$,~so
\[
-\varepsilon \mathsf{H}(x_{\varepsilon,u})\leq\underbrace{(c-\rB^{\intercal}u)^{\intercal }x_{\varepsilon,u}-\varepsilon\mathsf{H}(x_{\varepsilon,u})- (c-\rB^{\intercal}u)^{\intercal }x_{0,u}}_{=g_{\varepsilon}(u) - g_0(u)}+\varepsilon\mathsf{H}(x_{0,u})-\varepsilon\mathsf{H}(x_{0,u})\leq -\varepsilon \mathsf{H}(x_{0,u}).
\]
 It follows that, for any $u'\in\mathbb R^r$,  
\begin{equation}
\label{eq:boundGepsG0}
g_{\varepsilon}(u)-g_{\varepsilon}(u')\geq g_0(u)-g_0(u')-\varepsilon\left(\mathsf{H}(x_{\varepsilon,u})-\mathsf{H}(x_{0,u})\right).
\end{equation}
Instantiating this bound with any $u\in \mathcal C_{\text{subopt}}$ and $u'=u^{\star}$, where $\mathcal C_{\text{subopt}} = \cup_{j\in[V]\backslash \cI^{\star} }\mathcal C_j$ and $u^{\star}$ is some minimizer of $g_0$, we get
\[
g_{\varepsilon}(u)-g_{\varepsilon}(u^{\star})\geq g_0(u)-g_0(u^{\star})-\varepsilon R_{\mathsf H}
\]
so that if this lower bound is positive, we have $g_{\varepsilon}(u)>g_{\varepsilon}(u^{\star})$, and 
 $u\in \mathcal C_{\text{subopt}}$ cannot be optimal. To show this, recall from \eqref{eq:gCoercive} that $g_0$ is coercive so that it attains its minimum over $\mathcal C_{\text{subopt}}$ at some point $\bar u$ and, by definition of $\cI^{\star} $, $g_0(\bar u)-g_0(u^{\star})\geq \alpha$. Conclude that $\inf_{u\in }g_{\varepsilon}(u)>g_{\varepsilon}(u^{\star})\geq \inf_{u\in\mathbb R^r}g_{\varepsilon}(u)$ since $\varepsilon<\frac{\alpha}{R_{\mathsf H}}$, so that each global minimizer of $g_{\varepsilon}$ is an element of $\mathcal C_{\text{opt}}$, where $\cC_{\text{opt}} =\cup_{i\in \cI^{\star} } \mathcal C_i$. Recalling from \cref{prop:variationalFormND} that $\rB x_{\varepsilon}^{\star}$ minimizes $g_{\varepsilon}$, we have that $\rB x_{\varepsilon}^{\star}\in \cC_{\text{opt}}$.

We now control $\dist(\rB x_{\varepsilon}^{\star},\bd(\mathcal C_i))$. Recalling \eqref{eq:boundGepsG0}, it suffices to show that, for any $u\in \mathcal C_i$ satisfying $\dist(u,\bd(\mathcal C_i))<\dist(\rB v^{(i)},\bd(\mathcal C_i))/2=\frac{\gamma_i}{2}$, it holds that $g_0(u)-g_0(\rB v^{(i)})-\varepsilon R_{\mathsf H}>0$. Consequently, $g_{\varepsilon}(u)>g_\varepsilon(\rB v^{(i)})$ and thus $g_{\varepsilon}$ cannot attain its minimum at a point that is a distance less than $\dist(\rB v^{(i)},\bd(\mathcal C_i))/2$ from the boundary. For any such $u$,
\[
g_0(u)-g_0(\rB v^{(i)}) =\frac{1}{2} \|u\|^2 - u^\intercal \rB v^{(i)} - \frac{1}{2} \|\rB v^{(i)}\|^2 + \|\rB v^{(i)}\|^2= \frac{1}{2}\|u - \rB v^{(i)}\|^2.
\]
Now, if $\tilde u\in \bd(\mathcal C_i)$ is such that $\dist(u,\bd(\mathcal C_i))=\|\tilde u-u\|$, 
\[
\|u - \rB v^{(i)}\|\geq \|\tilde u-\rB v^{(i)}\|-\|\tilde u-u\|\geq \dist(\rB v^{(i)},\bd(\mathcal C_i))-\dist(u,\bd(\mathcal C_i))>\frac{\gamma_i}{2}.
\]
 Conclude that 
\[
g_{\varepsilon}(u)-g_{\varepsilon}(\rB v^{(i)})\geq g_0(u)-g_0(\rB v^{(i)})-\varepsilon R_{\mathsf H}>\frac{\gamma_i^2}{8}-\varepsilon R_{\mathsf H}>0,
\]
as $\varepsilon<\frac{\gamma_i^2}{8 R_{\mathsf H}}$. In sum,  $\dist(\rB x^{\star}_{\varepsilon},\bd(\mathcal C_i))\geq \dist(\rB v^{(i)},\bd(\mathcal C_i))/2$ so that $\rB x_{\varepsilon}^{\star}\in\Int(\mathcal C_i)$.

This bound will further enables us to control the constant from \cref{thm:entropicRateLP},
\[
    \kappa_{\left(c-\mathrm{B}^{\intercal}\mathrm{B}x_{\varepsilon}^{\star}\right)} = \min_{x\in\mathcal V(\mathcal K)\backslash \mathcal O_{c - \rB^\intercal \rB x^{\star}_{\varepsilon}}}(c-\mathrm{B}^{\intercal}\rB x^{\star}_{\varepsilon})^{\intercal}x-\min_{x\in\mathcal  K}(c-\mathrm{B}^{\intercal}\rB x^{\star}_{\varepsilon})^{\intercal}x.
\]
As $\rB x^{\star}_{\varepsilon}\in \Int(\mathcal C_i)$, $v^{(i)}\in \mathcal O_{c - \rB^\intercal \rB x^{\star}_{\varepsilon}}$. Now, fix  $v^{(j)}\in \argmin_{x\in\mathcal V(\mathcal K)\backslash \mathcal O_{c - \rB^\intercal \rB x^{\star}_{\varepsilon}}}(c-\mathrm{B}^{\intercal}\rB x^{\star}_{\varepsilon})^{\intercal}x$. Suppose, first, that $\mathrm{B}v^{(i)}\neq \mathrm{B}v^{(j)}$, then
\begin{equation}
\label{eq:distanceToBoundary}
\begin{aligned}
\dist(\rB x^{\star}_{\varepsilon},\bd(\mathrm{\mathcal C_i}))&\leq \dist(\rB x^{\star}_{\varepsilon},\{u\in\mathbb R^r: (c-\rB^{\intercal}u)^{\intercal}v^{(i)}=(c-\rB^{\intercal}u)^{\intercal}v^{(j)})\})
\\
&=\frac{|\rB(v^{(j)}-v^{(i)})^{\intercal} \rB x^{\star}_{\varepsilon}+c^{\intercal}(v^{(i)}-v^{(j)})|}{\|\rB(v^{(j)}-v^{(i)})\|}= \frac{\kappa_{\left(c-\mathrm{B}^{\intercal}\mathrm{B}x_{\varepsilon}^{\star}\right)}}{\|\rB(v^{(j)}-v^{(i)})\|},
\end{aligned}
\end{equation}
where the inequality follows by noting that the closest boundary point in $\cC_i$ to $\rB x_{\varepsilon}^{\star}$ is either an element of $\cC_i\cap \cC_j$ in which case equality holds, or it is not and the inequality is strict %
(recalling \cref{prop:minimizersG0}).
The first equality then follows by computing the projection onto a hyperplane.
Conclude that 
\[
\kappa_{\left(c-\mathrm{B}^{\intercal}\mathrm{B}x_{\varepsilon}^{\star}\right)}\geq \frac{\gamma_i}{2}\min_{\substack{k \in [V] \\ \rB(v^{(k)} - v^{(i)}) \neq 0}} 
\left\| \rB v^{(k)} - \rB v^{(i)} \right\|
\]in this case.
On the other hand, if $\rB v^{(i)}= \rB v^{(j)}$, $\kappa_{\left(c-\mathrm{B}^{\intercal}\mathrm{B}x_{\varepsilon}^{\star}\right)} = c^{\intercal} (v^{(j)}-v^{(i)})=\kappa_{\left(c-\mathrm{B}^{\intercal}\mathrm{B}x_{\varepsilon}^{\star}\right)}$ so that $\kappa_{\left(c-\mathrm{B}^{\intercal}\mathrm{B}x_{\varepsilon}^{\star}\right)}\geq \Delta_i$.

With these preparations, the claimed result follows by noting that 
   \begin{align*}
         \dist\left(x^\star_\varepsilon,\argmin_{x \in \cK}\left\{-\frac{1}{2}x^\intercal\rM x + c^\intercal x\right\}\right) &\leq  \dist\left(x^\star_\varepsilon,\argmin_{x \in \cK}\left(c-\rB^{\intercal}\rB v^{(i)}\right)^{\intercal}x\right)\\
         & = \dist\left(x^\star_\varepsilon,\argmin_{x \in \cK}(c-\rB^{\intercal}\rB x^{\star}_{\varepsilon})^{\intercal}x\right)\\
         &\leq 2R_1 \exp\left(\frac{-\kappa_{\left(c-\mathrm{B}^{\intercal}\mathrm{B}x_{\varepsilon}^{\star}\right)}}{\varepsilon R_1}+\frac{R_1 + R_H}{R_H}\right)\\
         &\leq 2 R_1 \exp\left(
  \frac{-\Delta_i}{\varepsilon R_1}
  + \frac{R_1 + R_{\sH}}{R_1}
\right),
    \end{align*}
   where the first inequality follows from \cref{prop:variationalFormND}, and the equality is due to the fact that $\rB v^{(i)},\rB x^\star_\varepsilon\in\Int (\mathcal C_i).$ The penultimate inequality follows by applying \cref{thm:entropicRateLP}. Finally, to obtain a uniform bound it suffices to take the minimum over all $i\in \cI^{\star} $ in the above bound and in the constraint on $\varepsilon$, which thus cannot be larger than $\min_{i\in I^\star}\frac{\gamma_i^2}{8 R_{\mathsf H}}$.
\qed

\subsection{Proofs for \texorpdfstring{\cref{sec:cert}}{Section 4}}
\subsubsection{Proof of \texorpdfstring{\cref{prop:nearCriticalPoints}}{Proposition 3}}
\label{proof:prop:nearCriticalPoints}
The result is based on the following pointwise estimates of the gradients. Its proof is included in \cref{proof:lem:gradientBound}.

\begin{lemma}
\label{lem:gradientBound}
    Fix $u\in\cC_i$, for some $i\in[V]$. Then, for any choice of $0<\delta<4\|\rB\|_{\mathrm{op}}R_1$ and 
    \[
0<\varepsilon\leq K_i\left( \dist(u,\bd(\cC_i))\right)\left( R_1+R_H-R_1\log\left(\frac{\delta}{4\|\rB\|_{\op}R_1}\right)\right)^{-1},
    \] 
    we have that
$
      \|\nabla g_{\varepsilon}(u)-\nabla g_{0}(u)\| \leq \delta/2$.
\end{lemma}

    As $K_i$ is an increasing function on $(0,\infty)$ and 
\[0<\varepsilon\leq K_i\left( \gamma\right)\left( R_1+R_H-R_1\log\left(\frac{\delta}{4\|\rB\|_{\op}R_1}\right)\right)^{-1},\]
the gradient estimate from \cref{lem:gradientBound} holds for every $u\in\cC_i$ with $\dist(u,\bd(\cC_i))\geq \gamma$.  
With this, any $u\in \cC_i$ satisfying $\| u-\rB v^{(i)}\|\leq \frac \delta 2$ and $\dist(u,\bd(\cC_i))\geq \gamma$ further has
    \[
        \|\nabla g_{\varepsilon}(u)\|\leq \|\nabla g_0(u)\|+\|\nabla g_{\varepsilon}(u)-\nabla g_{0}(u)\| \leq  \| u-\rB v^{(i)}\|+\frac{\delta}{2}\leq \delta,
    \]
   where we have used the fact that $\nabla g_{0}(u)=u-\rB v^{(i)}$ since $u\in \Int(\cC_i)$. 
  To prove the first assertion in Item (1), recall that $D_i = \dist(\rB v^{(i)},\bd(\cC_i))$ and observe that
\[
\dist(u,\bd(\cC_i))\geq D_i-\|u-\rB v^{(i)}\|\geq D_i-\frac{\delta}2
\]
by the reverse triangle inequality,
so that $\dist(u,\bd(\cC_i))\geq \gamma$ so long as  $\delta<D_i$ and $\gamma <D_i/2$. 

For the remainder of Item (1), note  that (a) and (b) together are equivalent to the statement that if $\rB v^{(i)}\in \cC_i,$ $\dist(u,\bd(\cC_i))\geq \gamma$, and $\|u-\rB v^{(i)}\|\geq \frac {3\delta} 2$, then, necessarily, $\|\nabla g_{\varepsilon}(u)\|>\delta$. This latter claim is a direct consequence of \cref{lem:gradientBound} and the reverse triangle inequality,  
    \begin{equation}
    \label{eq:boundGradients}
        \|\nabla g_{\varepsilon}(u)\|\geq \|\nabla g_0(u)\|-\|\nabla g_{\varepsilon}(u)-\nabla g_{0}(u)\| > \| u-\rB v^{(i)}\|-\frac{\delta}{2}\geq \delta.
    \end{equation}
    {Furthermore, if $\delta<D_i/3$, the points satisfying (a) are a subset of the open ball of radius $D_i/2$ so that $\dist(u,\bd(\cC_i))> D_i-D_i/2=D_i/2>\gamma$ and hence are disjoint from the points satisfying (b).}

    For Item (2), note that if $\rB v^{(i)}\not\in \Int(\cC_i)$ and $\dist(u,\bd(\cC_i))\geq\gamma$, then 
    \[
    \|u-\rB v^{(i)}\|\geq \dist(u,\bd(\cC_i))\geq \gamma \geq \frac{3}{2}\delta.
    \]
    Inserting the latter into \eqref{eq:boundGradients} yields $\|\nabla g_{\varepsilon}(u)\|>\delta$, concluding the proof.
    \qed

\subsubsection{Proof of \texorpdfstring{\cref{prop:alg_good}}{Proposition 4}}
\label{proof:prop:alg_good}
To prove that the algorithm terminates correctly with probability~$1$, we analyze the different termination conditions. We recall from \cref{prop:minimizersG0} that all local minimizers of $g_0$ must, necessarily, be interior points of $\cC_i$, for some $i\in[V]$. 

First, if $\bar u\in\Int(\cC_i)$ and $(c-\rB^{\intercal} \bar u)^{\intercal} \bar x > \min_{x\in \mathcal K} (c-\rB^{\intercal} \bar u)^{\intercal} x$, then we must have that $\rB \bar x \neq \rB v^{(i)}$, whence $\bar u$ is not a local minimizer of $g_0$ and \cref{alg:1} terminates at line $2$. Indeed, supposing to the contrary that $\rB \bar x = \rB v^{(i)}$ and that $(c-\rB^{\intercal} \bar u)^{\intercal} \bar x >  (c-\rB^{\intercal} \bar u)^{\intercal} v^{(i)}$, we obtain that  $c^{\intercal} \bar x>c^{\intercal} v^{(i)}$. As such, $(c-\rB^{\intercal} u_{\varepsilon})^{\intercal} \bar x >  (c-\rB^{\intercal} u_{\varepsilon})^{\intercal} v^{(i)}$, contradicting the fact that $\bar x \in \argmin_{x\in\mathcal K}(c-\rB^{\intercal}u_{\varepsilon})^{\intercal} x$.    

Assuming that the algorithm progresses past line $2$, we analyze the remaining termination conditions. Our approach is based on the following technical lemmas. The proofs of Lemmas \ref{lem:derivativesInterior}-\ref{lem:interiorIFF} are provided, respectively, in Appendices \ref{proof:lem:derivativesInterior}-\ref{proof:lem:interiorIFF}.
 
\begin{lemma}
    \label{lem:derivativesInterior}
   If $u\in \Int(\cC_i)$, $\rB \argmin_{x\in\mathcal K}(c-\rB^{\intercal}u)^{\intercal} x =\{\rB v^{(i)}\}$. 
\end{lemma}

\begin{lemma}
\label{lem:boundaryCondition1}
   Suppose that $u\in \cap_{r\in\mathcal R}\cC_r$ for some index set $\mathcal R\subset [V]$ with the property that $\cC_i\neq \cC_j$, for any $i,j\in\mathcal R$ with $i\neq j$. Then, for almost every $w$ sampled uniformly at random from the unit sphere, and every $\eta>0$, $u_+=u+\eta w$ and $u_-=u-\eta w$ are such that $u_+\in \cC_+,u_-\in\cC_-$ for distinct sets $\cC_+,\cC_-\in\{\cC_i\}_{i\in [V]}$, $u_+,u_-\not\in \cC_+\cap\cC_-$, and $u_+,u_-\not\in \cC_i\cap\cC_j$ for each $i,j\in\mathcal R$ with $i\neq j$.   
\end{lemma}

\begin{lemma}
\label{lem:boundaryPointLemma}
Suppose that $u',u''\in\mathbb R^r$ satisfy the property that 
$\rB x'\neq \rB x''$ and
\[
    (c-\rB^{\intercal} u)^{\intercal} x'= (c-\rB^{\intercal} u)^{\intercal} x'' = \min_{x\in\mathcal K}(c-\rB^{\intercal}u)^{\intercal} x, 
\]
where $x'\in\argmin_{x\in\mathcal K}(c-\rB^{\intercal}u')^{\intercal} x$ and $x''\in\argmin_{x\in\mathcal K}(c-\rB^{\intercal}u'')^{\intercal} x$. Then, $u\in\bd(\cC_i)$ for some  $i\in[V]$. 
\end{lemma}

\begin{lemma}
\label{lem:interiorIFF}
    We have $u\in\Int(\cC_i)$ for some $i\in[V]$ if and only if for almost every $w$ sampled uniformly from the unit sphere, $u_+=u+\eta w$ and $u_-=u-\eta w$ are such that $\rB x_+=\rB x_-$ for some $\eta>0$, where $x_+\in\argmin_{x\in\mathcal K}(c-\rB^{\intercal}u_+)^{\intercal} x$ and $x_-\in\argmin_{x\in\mathcal K}(c-\rB^{\intercal}u_-)^{\intercal} x$.     
\end{lemma}

With these lemmas in hand, we conclude the proof. If $\bar u\in \Int(\cC_i)$, \cref{lem:boundaryPointLemma}  implies that, for any given $\eta>0$, $u_+$ and $u_-$ are such that $\rB x_+=\rB x_-$ or $s_+=(c-\rB^{\intercal} \bar u)^{\intercal} x_+,s_-=(c-\rB^{\intercal} \bar u)^{\intercal} x_-,\text{ and } s=\min_{x\in\mathcal K}(c-\rB^{\intercal}\bar u)^{\intercal} x$ do not all take the same value. Conclude that the termination condition on line 14 ($s_+=s_-=s$) cannot be met in this setting unless, also,  $\rB x_+=\rB x_-$ which corresponds to the termination condition on line 10 (i.e., the algorithm  correctly returns that $\bar u$ is an interior point; recall \cref{lem:interiorIFF}). In light of \cref{lem:derivativesInterior}, once $x_+,x_-\in\Int(\cC_i)$, we have $\rB x_+=\rB x_-$. As $\|\bar u-u_+\|= \|\bar u-u_-\|=\frac{\eta_0}{2^{k-1}}$ at the  start of the $k$-th iteration of the while loop, $u_+$ and $u_-$ are interior points and the algorithm terminates after at most $\max\left\{0,\lceil \log_2(\eta_0/\zeta(\bar u))\rceil\right\}+1$ iterations. 

If $\bar u \in \bd(\cC_i)$, \cref{lem:interiorIFF} implies that the termination condition $\rB x_-=\rB x_+$ on line 10 is not satisfied with probability $1$. To see that $s_+=(c-\rB^{\intercal} \bar u)^{\intercal} x_+,s_-=(c-\rB^{\intercal} \bar u)^{\intercal} x_-,\text{ and } s=\min_{x\in\mathcal K}(c-\rB^{\intercal}\bar u)^{\intercal} x$ coincide  for some $\eta>0$, recall that at the start of the $k$-th iteration of the while loop, $\|\bar u-u_+\|= \|\bar u-u_-\|=\frac{\eta_0}{2^{k-1}}$ so that  after at most $\max\left\{0,\lceil \log_2(\eta_0/\lambda(\bar u))\rceil\right\}+2$ iterations, $\|\bar u-u_-\|=\|\bar u-u_+\|< \dist\left(\bar u, \cup_{\{j \in [V] : \bar u \not\in \mathcal C_j\}} \mathcal C_j\right)$. Letting $\mathcal R\coloneqq\{j \in [V] : \bar u \in \mathcal C_j\text{ and }\mathcal C_i\neq \mathcal C_j\text{ for every }i<j\}$, we obtain that   $u_-,u_+\in \Int(\cup_{r\in\mathcal R} \mathcal C_j)$ and, applying \cref{lem:boundaryCondition1}, $u_+,u_-\not \in \mathcal C_i\cap \mathcal C_j$ for each $i,j\in\mathcal R$ with $i\neq j$. In sum, $u_+$ and $u_-$ are interior points of $\cC_{r_+}$ and $\cC_{r_-}$, respectively, for some ${r_+,r_-\in\mathcal R}$ and $\bar u\in \cC_{r_+}\cap\cC_{r_-}$ so that 
\[
    (c-\rB^{\intercal}\bar u)^{\intercal} v^{(r_+)}=(c-\rB^{\intercal}\bar u)^{\intercal} v^{(r_-)}=\min_{x\in\mathcal K}(c-\rB^{\intercal}\bar u)^{\intercal} x.
\]
As $u_+\in\Int(\cC_{r_+})$ and $u_-\in\Int(\cC_{r_-})$, \cref{lem:derivativesInterior} implies that $\rB x_+=\rB v^{(r_+)}$, $\rB x_-=\rB v^{(r_-)}$ and, necessarily, $c^{\intercal} x_+=c^{\intercal} v^{(r_+)}$ and $c^{\intercal} x_-=c^{\intercal} v^{(r_-)}$ (otherwise $x_+$ and $x_-$ are suboptimal for the relevant problems, see line 9 of \cref{alg:1}). Inserting these equalities into the above display yields that $s_+=s_-=s$.  

 The final assertion in \cref{prop:alg_good} then follows by taking the maximum over all critical points, noting that there are only finitely many, as stated next (see \cref{proof:lem:finiteCardinality} for the proof).  
\begin{lemma}
\label{lem:finiteCardinality}
   The set of critical points of $g_0$ has finite cardinality. 
\end{lemma}
\qed

\subsection{Proof of \texorpdfstring{\cref{prop:slowRate}}{Proposition 5} }
\label{proof:prop:slowRate}
    Given that $x_{\varepsilon}^{\star}$ is optimal for the regularized problem,  we have that 
    \[
        \frac 12 \left(x^{\star}_{\varepsilon}\right)^{\intercal} \rM x^{\star}_{\varepsilon} + c^{\intercal} x^{\star}_{\varepsilon} - \varepsilon\mathsf H\left(x^{\star}_{\varepsilon}\right) \leq  \frac 12 \left(x^{\star}\right)^{\intercal} \rM x^{\star} + c^{\intercal} x^{\star} - \varepsilon\mathsf H\left(x^{\star}\right), 
    \]
    for any $x^{\star}\in \argmin_{x\in\mathcal K}\left\{\frac 12 x^{\intercal} \rM x + c^{\intercal} x\right\}$.
    Conclude that 
    \[
        \left(\frac 12 \left(x^{\star}_{\varepsilon}\right)^{\intercal} \rM x^{\star}_{\varepsilon} + c^{\intercal} x^{\star}_{\varepsilon}\right)- \left(\frac 12 \left(x^{\star}\right)^{\intercal} \rM x^{\star} + c^{\intercal} x^{\star}\right) \leq  \varepsilon \left(\mathsf H\left(x^{\star}_{\varepsilon}\right)-\mathsf H\left(x^{\star}\right)\right)\leq \varepsilon\sup_{x,x'\in\mathcal K}\left\{\mathsf H(x)-\mathsf H(x')\right\}, 
    \]
    which implies that the difference of costs is of order $O(\varepsilon)$. If, in addition, $x^{\star}_{\varepsilon}\to x^{\star}$, the first upper bound vanishes as $\varepsilon\to 0$, proving the second claim.  \qed

\section{Concluding Remarks}\label{sec: summary}

This work studied the approximation properties of entropically penalized QPs. We showed that concave QPs are unique in that the entropically regularized solutions and corresponding value of the quadratic cost converge to the solution set of the original problem at an exponential rate. By contrast, indefinite and concave QPs were shown to generally follow a much slower rate of $\Omega(\varepsilon)$ for the convergence of solutions. The exponential convergence rates where derived by leveraging a  new variational formulation of quadratic QPs, which tie them to a class of LPs with varying cost but fixed constraint set.

Building on this analysis, we showed that certain local solutions of the regularized variational problem are close to local minimizers of the unregularized variational objective. This observation gave rise to our \cref{alg:1}, which, under certain conditions, takes a $\delta$-critical point of the regularized variational objective and returns a local minimizer of its unregularized counterpart. Finally, these results were applied to the problem of GW alignment with kernel similarity measures. We recast the GW problem in a variational form compatible with our QP analysis, and establish, for the first time, an exponential decay rate for the entropic approximation gap between solution and alignment cost values. This is a notable improvement over the only previously known approximation bound of $\varepsilon\log(1/\varepsilon)$ \cite{zhang2024gromov}. %
Below, we outline some future research directions stemming from this work.

\medskip
\noindent\textbf{Efficient algorithms for the variational form.} As noted in \cref{rem:entropic_computation}, the variational form of the Euclidean EGW problem derived in \cite{zhang2024gromov} quickly inspired new algorithms for computing it subject to nonasymptotic convergence guarantees \cite{rioux2023entropic}. This algorithm further benefits from the fact that the regularized LPs can be solved faster than standard LPs via Sinkhorn's fixed-point iterations, which leverages the structure of the set of couplings of probability measures. It is of interest to identify other constraint sets which can benefit similarly from entropic regularization. Given such solvers, our variational form would enable efficient algorithms for optimizing general classes of concave QPs.

\medskip
\noindent{\textbf{Extending the variational form.}} Building on the initial results of \cite{zhang2024gromov}, we view \cref{prop:EGWDuality} as a second step towards a general variational theory for GW problems, effectively connecting them back to the well-studied OT problem. Already, the work of \cite{zhang2024gromov} has lead to significant statistical and computational progress for the quadratic GW distance.  Though the results derived herein are limited to the finite-dimensional case with a NSD cost matrix, we conjecture that similar representations can be derived for arbitrary cost matrices, still under the finitely discrete setting, using similar techniques. We are hopeful that these developments will provide useful insight into the infinite-dimensional setting as is the case for OT duality (see Section 1.2 in \cite{evans1997partial}). 

\medskip
\noindent\textbf{Improved rates for indefinite and convex QPs.}   As illustrated by the proof technique of \cref{thm:EntropicQP}, obtaining exponential convergence rates is contingent upon showing that all global solutions of the variational problem are interior points of some $\mathcal C_i$. As such, it appears reasonable to expect exponential rates for convex and indefinite QPs, under conditions enforcing this behavior once the relevant variational formulations have been established. Identifying such primitive conditions, i.e., which can be verified \emph{a priori}, is key for progress in this regard and forms an appealing research direction.      
\pagebreak

\bibliographystyle{amsalpha}
\bibliography{ref}

\appendix

\section{Proofs of  Lemmas}
\label{sec:proofsLemmas}
\subsection{Proof of \texorpdfstring{\cref{lem:gradientBound}}{Lemma 1}}
\label{proof:lem:gradientBound}
   First, observe that  the constant $\kappa_{c-\rB^{\intercal}u}$ from \cref{thm:entropicRateLP} satisfies $\kappa_{c-\rB^{\intercal}u}\geq K_i\left(\dist(u,\bd(\cC_i))\right)$
 which follows from \eqref{eq:distanceToBoundary} and the surrounding discussion. {Note that $K_i(\gamma)$ plays the role of  the constant $\Delta_i$ in \cref{thm:informalRate}. However, the constant $\kappa_{c-\rB^{\intercal}\rB v^{(i)}}$ from $\Delta_i$ must be replaced with the minimum of $c^{\intercal}(v^{(k)}-v^{(i)})$ over vertices $v^{(k)}$ satisfying $\rB v^{(k)}=\rB v^{(i)}$ and $c^{\intercal}v^{(k)}\neq c^{\intercal}v^{(i)}$ as $\rB v^{(i)}$ may not be an element of $\cC_i$.} Now, since $\delta<4\|\rB\|_{\mathrm{op}}R_1$, we obtain
\[
\varepsilon\leq K_i\left( \dist(u,\bd(\cC_i))\right)\left( R_1+R_H\right)^{-1}\leq \frac{\kappa_{c-\rB^{\intercal}u}}{R_1+R_{\sH}},  
\]
so that \cref{thm:entropicRateLP} can be applied, yielding \[
     \|\nabla g_{\varepsilon}(u)-\nabla g_{0}(u)\|=\|\rB x_{u,\varepsilon}-\rB v^{(i)}\| \leq 2\|\rB\|_{\mathrm{op}} R_1 \exp\left(-\frac{\kappa_{c-\rB^{\intercal}u}}{\varepsilon R_1}+\frac{R_1+R_{\mathsf{H}}}{R_1}\right),
     \]
     where $\{x_{u,\varepsilon}\}=\argmin_{x\in\mathcal K}\left\{ (c-\rB^{\intercal}u)^{\intercal}x-\varepsilon\sH(x)\right\}$.

The result then follows by noting that 
\[
\frac{\kappa_{c-\rB^{\intercal}u}}{\varepsilon}\geq  \frac{\kappa_{c-\rB^{\intercal}u}}{K_i\left( \dist(u,\bd(\cC_i))\right)}\left( R_1+R_H-R_1\log\left(\frac{\delta}{4\|\rB\|_{\op}R_1}\right)\right)\geq  R_1+R_H-R_1\log\left(\frac{\delta}{4\|\rB\|_{\op}R_1}\right),
\]
for
$0<\varepsilon\leq K_i\left( \dist(u,\bd(\cC_i))\right)\left( R_1+R_H-R_1\log\left(\frac{\delta}{4\|\rB\|_{\op}R_1}\right)\right)^{-1}$.
\qed

\subsection{Proof of \texorpdfstring{\cref{lem:derivativesInterior}}{Lemma 2} }
\label{proof:lem:derivativesInterior}
  As $u\in\Int(\cC_i)$, there exists a neighborhood of $u$ on which $g_0=\frac 12 \|\cdot\|^2+(c-\rB^{\intercal}(\cdot))^{\intercal} v^{(i)}$ so that $g_0$ is smooth at $u$ with derivative $u-\rB v^{(i)}$. On the other hand, as noted in \cref{rem:entropic_computation}  the subdifferential of $g_0$ at any point $u'$ is given by $u'-\rB \argmin_{x\in\mathcal K}(c-\rB^{\intercal}u')^{\intercal}x$. Applying Proposition 2.2.4 and the Corollary on p.32 in \cite{clarke1990optimization}, we obtain that $\rB \argmin_{x\in\mathcal K}(c-\rB^{\intercal}u)^{\intercal}x=\{\rB v^{(i)}\},$ proving the claim.   \qed
\subsection{Proof of \texorpdfstring{\cref{lem:boundaryCondition1}}{Lemma 3} }
\label{proof:lem:boundaryCondition1}
   Fix a unit vector $w$ and suppose that $w^{\intercal}\rB v^{(i)}<w^{\intercal}\rB v^{(j)}$ for some $i,j\in\mathcal R$. Then, 
   \[
        (c-\rB^{\intercal} u_+)^{\intercal} \rB v^{(i)}= (c-\rB^{\intercal} u)^{\intercal} \rB v^{(i)} +\eta w^{\intercal}\rB v^{(i)} = (c-\rB^{\intercal} u)^{\intercal} \rB v^{(j)} +\eta w^{\intercal}\rB v^{(i)}< (c-\rB^{\intercal} u_+)^{\intercal} \rB v^{(j)} 
   \]
   and, likewise, $(c-\rB^{\intercal} u_-)^{\intercal} \rB v^{(j)}<(c-\rB^{\intercal} u_-)^{\intercal} \rB v^{(i)}$. It follows that $u_+$ and $u_-$ must lie in distinct elements $\cC_+$ and $\cC_-$ of $\{\cC_i\}_{i\in [V]}$ and $u_+,u_-\not\in \cC_+\cap\cC_-$. 

   It remains to show that a uniform sample, $w$, from the unit sphere in $\mathbb R^r$ satisfies the property that $w^{\intercal}\rB v^{(i)}\neq w^{\intercal}\rB v^{(j)}$ for each $i,j\in\mathcal R$ with $i\neq j$ with probability $1$. It is well-known (cf. e.g., Problem 20.26 in \cite{billingsley2012probability}) that if $Z$ is distributed according to the standard Gaussian distribution in $\mathbb R^r$, then $Q = Z/\|Z\|$ is distributed according to the uniform measure on the unit sphere. Conclude that, for any $a\in\mathbb R^r\backslash \{0\}$, $\langle Q,a\rangle\neq 0$ with probability $1$ as $\langle Q,a\rangle$ is distributed according to a nondegenerate univariate Gaussian random variable. It follows that $w^{\intercal}\rB v^{(i)}\neq w^{\intercal}\rB v^{(j)}$ almost surely for each $i,j\in\mathcal R$ with $i\neq j$  as $\rB v^{(i)}\neq \rB v^{(j)}$ for each $i,j\in\mathcal R$ with $i\neq j$ as $\cC_i$ and $\cC_j$ are assumed to be distinct.\footnote{If $\rB v^{(i)}=\rB v^{(j)}$, either $c^{\intercal} v^{(i)}=c^{\intercal} v^{(j)}$ in which case $\cC_i=\cC_j$ or $c^{\intercal} v^{(i)}\neq c^{\intercal} v^{(j)}$ in which case $\cC_i$ or $\cC_j$ is empty.}      
\qed

\subsection{Proof of \texorpdfstring{\cref{lem:boundaryPointLemma}}{Lemma 4} }
\label{proof:lem:boundaryPointLemma}
Suppose to the contrary that $u$ is an interior point. As
$
    (c-\rB^{\intercal} u)^{\intercal} x'= (c-\rB^{\intercal} u)^{\intercal} x'' = \min_{x\in\mathcal K}(c-\rB^{\intercal}u)^{\intercal} x, 
$ we obtain that
$x',x''\in\argmin_{i\in[V]}(c-\rB^{\intercal}u)^{\intercal} x$. However, \cref{lem:derivativesInterior} then asserts that $\rB x'=\rB x''$ contradicting the assumption that $\rB x'\neq \rB x''$. 
\qed

\subsection{Proof of \texorpdfstring{\cref{lem:interiorIFF}}{Lemma 5}}
\label{proof:lem:interiorIFF}
    If  $u\in \bd(\cC_i)$, \cref{lem:boundaryCondition1} yields that $u_+$ and $u_-$ lie, respectively, in distinct sets $\cC_+,\cC_-\in\{\cC_i\}_{i\in[V]}$ and $u_+,u_-\not \in \cC_+\cap \cC_-$ for every $\eta>0$ with probability $1$. Now, note that if  $\rB x_+=\rB x_-$,
$
(c-\rB^{\intercal}u_+)^{\intercal} x_+= c^{\intercal} x_+-u_+^\intercal\rB x_-
$ and $c^{\intercal} x_-=c^{\intercal} x_+$ as if $c^{\intercal} x_-<c^{\intercal} x_+$, $x_+\not\in\argmin_{x\in\mathcal K}(c-\rB^{\intercal}u_+)^{\intercal} x$ and an analogous argument holds if $c^{\intercal} x_->c^{\intercal} x_+$. Conclude that $x_-\in\argmin_{x\in\mathcal K}(c-\rB^{\intercal}u_+)^{\intercal} x$ so that the problems $\min_{x\in\mathcal K}(c-\rB^{\intercal}u_+)^{\intercal} x$ and $\min_{x\in\mathcal K}(c-\rB^{\intercal}u_-)^{\intercal} x$ share at least one vertex solution. It follows that    either $u_+,u_-\in\cC_+\cap\cC_-$ or $\cC_+=\cC_-$, contradicting \cref{lem:boundaryCondition1}. On the other hand, if $u\in\Int(\cC_i)$, then, for every $\eta>0$ sufficiently small, $u_+,u_-\in\Int(\cC_i)$ so that \cref{lem:derivativesInterior} implies that $\rB x_+=\rB x_-$.  
\qed
\subsection{Proof of  \texorpdfstring{\cref{lem:finiteCardinality}}{Lemma 6}}
 \label{proof:lem:finiteCardinality} 
  Recall from \eqref{eq:subdifferential} that $u$ is a critical point of $g_0$ if and only if $ u =\rB x_{u}$ for some $x_{u}\in\argmin_{x\in\mathcal K}(c-\rB^{\intercal} u)^{\intercal} x$. If $ u$ is an interior point of some $\cC_i$ for $i\in[V]$, we have from \cref{lem:derivativesInterior} that $\rB x_{ u}=\rB v^{(i)}$ so that $g_0$ admits at most one critical point in $\Int(\cC_i)$, namely $\rB v^{(i)}$.

   Now, if $u$ is both a critical point and a boundary point, let $I(u) = \{i\in[V]:u \in \mathcal C_i\}$ which contains at least two elements.
    We now show that $u$ is the unique critical point with this value of $I(u)$. Indeed, suppose to the contrary that $I(u+v)=I(u)$ for some $v\in\mathbb R^r$. Then, for every $i,j\in I(u)$ it holds that 
    \[
        (c-\mathrm{B}^{\intercal}  u)^{\intercal} x^{(i)} = (c-\mathrm{B}^{\intercal}  u)^{\intercal} x^{(j)},  
   \text{  
    and } 
        (c-\mathrm{B}^{\intercal} ( u+v))^{\intercal} x^{(i)} = (c-\mathrm{B}^{\intercal} ( u+v))^{\intercal} x^{(j)}
        \]
        so that
         $v^{\intercal} \mathrm Bx^{(i)} = v^{\intercal} \mathrm Bx^{(j)}\eqqcolon \vartheta .$
    By virtue of being critical points, there exists nonnegative $(\lambda_i)_{i\in I(u)}, (\lambda_i')_{i\in I(u)}$ summing to one for which 
    \[
         u =\sum_{i\in I(u)} \lambda_i\mathrm{B}x^{(i)},\quad  u+v =\sum_{i\in I(u)} \lambda_i'\mathrm{B}x^{(i)}\text{ whereby } v = \sum_{i\in I(u)}(\lambda_i'-\lambda_i)\mathrm{B}x^{(i)}. 
    \]
    Taking the inner product of the rightmost expression with $v$ yields that 
    \[
        \|v\|^2 = \sum_{i\in I(u)}(\lambda_i'-\lambda_i)v^{\intercal}\mathrm{B}x^{(i)} = \vartheta \sum_{i\in I(u)}(\lambda_i'-\lambda_i)=0, 
    \]
    so that $v$ must be $0$, proving the claimed result.
\qed

\section{Primer on Optimal Transport and Gromov-Wasserstein Problems}
\label{sec:primerOT}
\subsection{Optimal Transport}
Optimal transport (OT) theory \cite{villani2003topics,villani2008optimal,santambrogio15} provides a general framework for comparing probability distributions. Given  two Borel probability measures $\mu_0,\mu_1$ on $\RR^d$ and a cost function $c:\RR^d\times \RR^d\to \RR$, the OT problem is defined by
\begin{equation}
    \mathsf{OT}_c(\mu_0,\mu_1)\coloneqq\inf_{\pi\in\Pi(\mu_0,\mu_1)}\int c(x,y)d\pi(x,y),\label{EQ:Kantorovich_OT}
\end{equation}
where $\Pi(\mu_0,\mu_1)$ is the set of couplings between $\mu_0$ and $\mu_1$. One of the most important instances of OT occurs when $c(x,y)=\|x-y\|^p$ for some $p\geq 1$. In this case $\mathsf{OT}_c(\mu_0,\mu_1)$ is the $p$-th power of the $p$-Wasserstein distance between $\mu_0$ and $\mu_1$. The $p$-Wasserstein distance defines a \emph{bona fide} metric on the space of probability measures with finite $p$-th moment and endows this space with a rich geometry, which has lead to the wide adoption of OT-based tools in machine learning \cite{arjovsky2017wasserstein,blanchet2019quantifying,courty2016optimal,chen2022inferential,gulrajani2017improved,rubner2000earth,tolstikhin2017wasserstein} and statistics \cite{bernton2019approximate,bernton2019parameter,bigot2020statistical,
backhoff2022bayesian,carlier2017vector,chernozhukov2017monge,chen2023wasserstein,ghosal2022multivariate,hallin2021quantile,
panaretos2020invitation,torous2021optimal,zhang2022wasserstein}.  

When $\mu_0$ and $\mu_1$ are supported on finite sets of points $\mathcal X_0$ and $\mathcal X_1$ with $|\mathcal X_0|=N_0$ and $|\mathcal X_1|=N_1$  respectively, \eqref{EQ:Kantorovich_OT} is  an LP in $N_0N_1$ variables and can be solved using the network simplex method, which has  a computational complexity of $O((N_0+N_1)N_0N_1)$ up to logarithmic factors \cite{orlin1997polynomial,tarjan1997dynamic}, see Section 3.5 in \cite{peyre2019computational} for details. To accelerate computation of OT it was proposed in \cite{cuturi2013sinkhorn} to regularize \eqref{EQ:Kantorovich_OT}, yielding the EOT problem \cite{leonard2012schrodinger,schrodinger1931uber}
\begin{equation}
    \mathsf{EOT}^{\varepsilon}_c(\mu_0,\mu_1):=\inf_{\pi\in\Pi(\mu_0,\mu_1)}\left\{\int c(x,y) d\pi(x,y) - \varepsilon\mathsf H(\pi)\right\}\label{EQ:EOT}
\end{equation}
where $\varepsilon>0$ is a regularization parameter and
\[
    \mathsf H(\pi) =- \sum_{(x_0,x_1)\in\mathcal X_0\times \mathcal X_1} \pi\left(\left\{\left(x_0,x_1\right)\right\}\right)\log\left(\pi\left(\left\{\left(x_0,x_1\right)\right\}\right)\right).
\] 
With this, \eqref{EQ:EOT} can be identified with an entropically penalized LP, \eqref{eq:entropicLP}, and can be solved using Sinkhorn's algorithm \cite{sinkhorn1967diagonal}, which has a per iteration computational complexity of $O(N_0N_1)$ (see Section 4.3 in \cite{peyre2019computational}); bounds on the total number of iterations required for convergence are provided in \cite{franklin1989scaling,knight2008sinkhorn}.   

\subsection{Gromov-Wasserstein Problems} While the standard OT formulation \eqref{EQ:Kantorovich_OT} can be extended to the case where $\mu_0$ and $\mu_1$ are measures on arbitrary Polish spaces, comparing distributions across different spaces requires defining a cost function on the product space. However, constructing costs which capture desirable notions of similarity between two spaces may be highly nontrivial. By contrast, GW problems \cite{chowdhury2019gromov,Memoli11,sturm2012space} enable comparing distributions given only similarity measures on each space (e.g., metrics) and serve as a quadratic counterpart to the linear OT problem. Namely, given probability measures $\mu_0,\mu_1$ on the Polish spaces $\mathcal X_0,\mathcal X_1$ and similarity measures $k_0:\mathcal X_0\times\mathcal X_0\to\mathbb R,k_1:\mathcal X_1\times\mathcal X_1\to\mathbb R$, the $(p,q)$-GW problem for $p,q>0$ is given by   

\begin{equation}
\label{eq:GromovWassersteinDefinition}
\mathsf{GW}_{p,q}(\mu_0,\mu_1)\coloneqq \left(\inf_{\pi\in\Pi(\mu_0,\mu_1)}\int| k_0^q(x,x')-k_1^q(y,y')|^pd\pi(x,y) \pi(x',y')\right)^{1/p}.
\end{equation}
This problem can be thought of, equivalently, as comparing the objects $(\mathcal X_0,k_0,\mu_0)$ and $(\mathcal X_1,k_1,\mu_1)$, which are known as metric measure (mm) spaces when $k_0$ and $k_1$ are metrics on the relevant spaces. Remarkably, $\mathsf{GW}_{p,q}$ defines a metric on the space of all mm spaces with finite $pq$-size (i.e., $\int k_0^{pq}d\mu_0d \mu_0,\int k_1^{pq}d\mu_1d\mu_1<\infty$) modulo the equivalence defined by isomorphism.\footnote{Two mm spaces $(\mathcal X_0,k_0,\mu_0)$ and  $(\mathcal X_1,k_1,\mu_1)$ are said to be isomorphic if there exists an isometry $T:\supp(\mu_0)\to \supp(\mu_1)$ for which $\mu_0\circ T^{-1}=\mu_1$.}
This distance serves as an OT-based $L^p$ relaxation of the classical Gromov-Hausdorff distance and affords both a means to compare mm spaces and to align them via a solution $\pi^{\star}$ to \eqref{eq:GromovWassersteinDefinition}. Thanks to these favorable properties, GW problems have recently been applied to numerous alignment tasks including single-cell genomics \cite{blumberg2020mrec,cao2022manifold,
demetci2020gromov}, alignment of language models \cite{alvarez2018gromov}, shape matching \cite{koehl2023computing,memoli2009spectral}, graph matching \cite{chen2020graph,petric2019got,xu2019scalable,xu2019gromov}, heterogeneous domain adaptation \cite{sejourne2021unbalanced,yan2018semi}, and generative modeling \cite{bunne2019learning}. 

When $\mu_0$ and $\mu_1$ are 
supported on finite sets of points $\mathcal X_0$ and $\mathcal X_1$ with $|\mathcal X_0|=N_0$ and $|\mathcal X_1|=N_1$  respectively
the GW problem is simply a QP, though it is generally nonconvex. As such, it belongs to the class of NP-hard problems \cite{pardalos1991quadratic}. To speed up computation,  \cite{peyre2016gromov,solomon2016entropic} propose to regularize  the GW problem as
\[
\mathsf{EGW}_{p,q}^\varepsilon(\mu_0,\mu_1)\coloneqq\left\{\inf_{\pi\in\Pi(\mu_0,\mu_1)} \iint \left|k_0^q(x,x')-k_1^q(y,y')\right|^p d\pi(x,y) \pi(x',y')+ \varepsilon {\mathsf {H}}(\pi)\right\},
\]
and propose an optimization procedure based on linearizing the quadratic term, so that each iteration of the method requires solving a single entropic OT problem. However, these works did not establish if the proposed algorithms terminate after a finite number of steps. 

When $\mathcal X_0=\mathbb R^{d_0},\mathcal X_1=\mathbb R^{d_1}$ and $\kappa_0$ and $\kappa_1$ are either the Euclidean norm or the inner product on the relevant spaces, \cite{rioux2023entropic} provided the first algorithms for solving the entropic $(2,2)$-GW problem with convergence rate results and the same per iteration complexity as a refined version of the method in \cite{peyre2016gromov} which exploits the low-rank structure of the cost (see \cite{scetbon2022linear} for details). This new algorithm consists of optimizing an alternative representation of the entropic GW problem, derived in \cite{zhang2024gromov},
\[
\mathsf{EGW}_{2,2}^{\varepsilon}(\mu_0,\mu_1)=
\mathsf K({\mu_0,\mu_1})+\mathsf R^{\varepsilon}(\mu_0,\mu_1),
\]
 where, for $\kappa_i(x,x')=\|x-x'\|$ for $i=0,1$,
\[
\begin{aligned}
&\mathsf K({\mu_0,\mu_1})=\iint\kappa^4_0d\mu_0d\mu_0+\iint\kappa_1^4d\mu_1d\mu_1-4\iint\|x\|^2\|y\|^2d\mu_0(x)d\mu_1(y),\\
      &\mathsf{R}^{\varepsilon}(\mu_0,\mu_1) = \inf_{\mathrm Z\in\mathbb R^{d_x\times d_y}}\left\{32\|\mathrm Z\|_{\mathrm F}^2+\inf_{\pi\in\Pi(\mu_0,\mu_1)}\int -4\|x\|^2\|y\|^2-32x^{\intercal}\mathrm Zyd\pi(x,y)+\varepsilon{\mathsf{H}}(\pi)\right\},    
\end{aligned}
\]
   for any $\varepsilon\geq 0$ and $\|\cdot\|_{\mathrm F}$ is the Frobenius norm. With this, it suffices to solve the optimization problem defining $\mathsf R^{\varepsilon}(\mu_0,\mu_1)$. A similar representation was derived in Appendix E of \cite{rioux2023entropic} for the inner product similarity measure. Remarkably, these variational formulations connect the GW problem to a family of OT problems, enabling us to leverage the tools developed in OT theory. Such variational formulations are generally unknown in other settings, owing to the nonconvex quadratic nature of the GW problem. 

\end{document}